\documentclass[12pt]{article}

\usepackage[margin=1in]{geometry}
\usepackage[T1]{fontenc}
\usepackage{lmodern}
\usepackage{titlesec}
\usepackage{setspace}

\usepackage{amsmath,amssymb,amsthm}
\newtheorem*{theorem*}{Theorem}
\usepackage{mathtools}
\usepackage{mathrsfs}

\usepackage{enumitem}
\usepackage{booktabs}
\usepackage{array}

\usepackage{xcolor}
\usepackage{hyperref}

\hypersetup{
    colorlinks=true,
    linkcolor=blue,
    citecolor=blue,
    urlcolor=blue,
    pdftitle={Local Least Periods and the Distribution of Periodic Points in Dynamical Systems},
    pdfauthor={Diyath Nelaka Pannipitiya},
    pdfsubject={Topological Dynamics and Periodic Points},
    pdfcreator={LaTeX}
}

\newtheorem{theorem}{Theorem}[section]
\newtheorem{lemma}[theorem]{Lemma}
\newtheorem{proposition}[theorem]{Proposition}
\newtheorem*{proposition*}{Proposition}
\newtheorem{corollary}[theorem]{Corollary}

\theoremstyle{definition}
\newtheorem{definition}[theorem]{Definition}
\newtheorem{example}[theorem]{Example}

\theoremstyle{remark}
\newtheorem{remark}[theorem]{Remark}
\newtheorem{question}[theorem]{Question}

\newcommand{\Per}{\operatorname{Per}}

\begin{document}

\title{Local Periodic Orbit Penetration in Dynamical Systems}

\author{Diyath Nelaka Pannipitiya}

\date{}

\maketitle

\begin{abstract}

We introduce two complementary frameworks for studying local periodicity and recurrence in dynamical systems. In the topological setting, we consider continuous maps with nonempty periodic-point sets, restricted to the closures of their periodic points. We define the periodic penetration time of an open set as the least period of a periodic point it contains. Using dynamically refined finite generating covers, we construct asymptotic mean and variance invariants that describe local periodic penetration relative to the logarithms of minimum subcover cardinalities. We establish their invariance under topological conjugacy and investigate their behavior for finite periodic systems, full shifts, shifts of finite type, and zero-entropy systems with dense periodic points. We also extend the construction to invertible systems using symmetric dynamical refinements.

In the measure theoretic setting, we define local return times through positive measure intersections of measurable sets with their iterated preimages. Using dynamically refined finite generating partitions, we construct corresponding asymptotic mean and variance invariants and establish their invariance under measure-theoretic isomorphism. We derive explicit formulas for finite probability-preserving systems and study finite periodic approximations of uniform Bernoulli shifts on $q\geq 2$ symbols, proving that their system-level mean and variance invariants converge to $1/\log q$ and $0$, respectively. Our results demonstrate that local periodic penetration and return-time statistics provide information about dynamical systems beyond that captured by classical entropy.
\end{abstract}



\section{Introduction}
\label{Section: Introduction}

Periodic points play a central role in the study of dynamical systems. Their existence, distribution, and periods provide information about the orbit structure of a system and its long term behavior. In particular, periodic points are closely
related to several fundamental notions in topological dynamics,
including expansivity, chaos theory (especially in Devaney's definition of chaos), symbolic dynamics, and topological
entropy.

For a continuous map $f:X\to X$, the density of periodic points
ensures that every nonempty open subset of $X$ contains a
periodic point. However, this property alone does not describe
the periods of the periodic points encountered in increasingly
small neighborhoods. In particular, it does not distinguish
between regions containing periodic points of relatively small
period and those in which only periodic points of large period
can be found.

This observation motivates the study of the least period
required for a periodic orbit to intersect a given open set.
Rather than counting periodic points or periodic orbits, we
investigate how these least periods behave under successive
dynamical refinements of finite generating covers. Our aim is
to construct quantities that describe both the average size
and the variation of local periodic penetration times.

A natural point of comparison is topological entropy, which
measures the exponential growth of the minimum cardinalities
of dynamically refined open covers. Our constructions use
these covering numbers as a normalization, while incorporating
additional information about periodic orbit structure.
Consequently, they can distinguish certain dynamical systems
having the same topological entropy.

We also develop a measure theoretic counterpart. In a
probability preserving dynamical system, individual periodic
points may have measure zero or may not exist at all. Therefore, unlike in the topological version of the theory, we do not need the condition that the set of periodic points of the system being non-empty. Accordingly, we replace periodic penetration times with positive measure overlap times, which describe the least
positive iterate for which a measurable set intersects its
preimage in positive measure. This leads to mean and variance
quantities associated with finite measurable partitions and
their dynamical refinements.

\subsection{Main Results}

We begin with the topological setting. Let $X$ be a nonempty
compact Hausdorff space and let $f:X\to X$ be continuous. By $\Per(f)$, let us denote the set of all periodic points of $f$ in $X$. We assume that $\displaystyle \Per(f)\neq \varnothing$ and work on $\overline{\Per(f)}$, so that periodic points are dense
in the space under consideration.

For each nonempty open set $U\subseteq X$, we define its
periodic penetration time by
\[\eta_f(U) = \min\{\operatorname{per}_f(x): x\in U\cap\Per(f)\}.\]
Using finite generating covers that satisfy a linear growth
condition on maximal penetration times, we introduce the
system-level asymptotic mean invariant $\Lambda_A(f)$ in Definition \ref{Definition of Lambda_A} and its variance counterpart $\Lambda_V(f)$ in Definition \ref{Definition Penetration-Variance}.

\noindent We establish that both quantities are invariant under
topological conjugacy:

\begin{theorem*}[Theorem  \ref{Theorem: Lambda is Conjugacy Invariant}]
Suppose that $\displaystyle h:(X,f)\longrightarrow(Y,g)$ is a topological conjugacy, and suppose both $\mathscr{A}_f$ and $\mathscr{A}_g$ (see Definition \ref{Definition: eta admissible generating covers}) are non empty. Then,
\[\Lambda_A(f)=\Lambda_A(g).\]    
\end{theorem*}

\begin{proposition*}[Proposition \ref{proposition variance conjugacy}]
Let $X$ and $Y$ be nonempty compact Hausdorff spaces, and let $f:X\to X$ and $g:Y\to Y$ be continuous maps with dense periodic points. Suppose that $(X,f)$ and $(Y,g)$ are topologically
conjugate and that $\mathscr A_f\neq\varnothing$ (see Definition \ref{Definition: eta admissible generating covers}). Then,
\[\Lambda_V(f)=\Lambda_V(g).\]
\end{proposition*}

\noindent For finite periodic systems, we obtain explicit formulas in terms of the lengths of their periodic orbits. These formulas show that the mean and variance invariants can detect differences in periodic orbit structure even when topological entropy vanishes.
\begin{proposition*}[Proposition \ref{calculating Lambda for finite systmes}]
Let $X$ be a finite dynamical system consisting of $r$ distinct
periodic orbits with least periods $m_1,\ldots,m_r$, and let
$f:X\to X$ be the corresponding map. Then,
\[\Lambda_A(f) = \frac{\displaystyle\sum_{i=1}^{r}m_i^2}
{\displaystyle \left(\sum_{i=1}^{r}m_i\right) \left(1+\log\left(\sum_{i=1}^{r}m_i\right)\right)}.\]
\end{proposition*}

\begin{proposition*}[Proposition \ref{proposition variance-finite}]
Suppose $X$ is a finite discrete space and $f:X\to X$ is a permutation whose distinct cycles have lengths $m_1,\ldots,m_r$. Set $\displaystyle M=\sum_{i=1}^r m_i$. Then,
\[\displaystyle  \Lambda_V(f)=  \frac{1}{(1+\log M)^2}
 \left(\frac{\displaystyle \sum_{i=1}^r m_i^3}{M} -\left(\frac{\displaystyle \sum_{i=1}^r m_i^2}{M}\right)^2 \right).\]
In particular, $\Lambda_V(f)=0$ if and only if all the cycle lengths $m_1,\ldots,m_r$ are equal.
\end{proposition*}

\noindent For the one-sided full shift on $q\geq2$ symbols, we show that the standard cylinder generating partition satisfies
\[L_A(\sigma_q,\mathcal C_1) = \frac{1}{\log q}\ \text{(Theorem \ref{L_A for One-Sided Full Shifts})} \qquad \text{and} \qquad L_V(\sigma_q,\mathcal C_1)=0.\ \text{(Proposition \ref{Proposition: variance fullshift})} \]

\noindent We also establish linear bounds on periodic penetration
times for mixing shifts of finite type (Subsection \ref{subsection: Mixing Shifts of Finite Type}) and investigate the behavior of the invariants in zero entropy systems with dense periodic points. In particular, we construct a zero-entropy example for which $\displaystyle \Lambda_A(f)=\infty$. (Subsection \ref{subsection: A Zero-Entropy Example with Dense Periodic Points}.)

We extend the topological construction to invertible dynamical systems by replacing forward dynamical joins with symmetric joins (Definition \ref{Symmetric Iterated Join}) involving both positive and negative iterates. This leads to corresponding mean and variance invariants for homeomorphisms, including the two-sided full shift:

\begin{theorem*}[Theorem \ref{Theorem: L_A(sigma,C_0)=1/log q}: Two-Sided Full Shifts]
For the two-sided full shift $(\Sigma_q,\sigma_q)$, the partition
$\mathcal C_0$ is a symmetric generating partition and is
symmetrically $\eta$-admissible. Moreover,
\[M_{\sigma_q}\left(\widehat{(\mathcal C_0)}_{\sigma_q}^{(n)}\right) = 2n+1 \qquad \text{and} \qquad \displaystyle \lim_{n \rightarrow \infty }\frac{A_{\sigma_q}
\left(\widehat{(\mathcal C_0)}_{\sigma_q}^{(n)}\right)}{2n+1} = 1.\]
Consequently, $\displaystyle \widehat{L}_A(\sigma_q,\mathcal C_0) = \frac{1}{\log q}$, and therefore $\displaystyle \widehat{\Lambda}_A(\sigma_q) \geq \frac{1}{\log q}$.
\end{theorem*}

\noindent In the measure-theoretic setting, let
$(X,\mathcal B,\mu,f)$ be a probability-preserving
dynamical system. For a measurable set $B$ of positive
measure, we consider
\[\tau_\mu(B) = \min\{k\geq1: \mu(B\cap f^{-k}B)>0\}.\]
Using dynamically refined finite measurable partitions,
we introduce asymptotic mean and variance invariants
$\Lambda_\mu(f)$ and $\Lambda_{V,\mu}(f)$, subject to
an appropriate return-time admissibility condition.

\noindent We prove that these invariants are preserved under measure-theoretic isomorphism and derive explicit formulas for finite probability preserving systems:
\begin{theorem*}[Theorem \ref{Theorem: Measure-Theoretic Isomorphism Invariance}: Measure-Theoretic Isomorphism Invariance]
Let $(X,\mathcal B,\mu,f)$ and $(Y,\mathcal C,\nu,g)$ be measure-theoretically isomorphic probability preserving dynamical systems.
Then, $(X,\mathcal B,\mu,f)$ admits an admissible finite generator (see Definition \ref{def: Return-Time Admissibility}) if and only if $(Y,\mathcal C,\nu,g)$ does. Furthermore, whenever the global invariants are defined, we have
\[\Lambda_\mu(f)=\Lambda_\nu(g) \qquad\text{and}\qquad \Lambda_{V,\mu}(f)=\Lambda_{V,\nu}(g).\]
\end{theorem*}

\noindent We also study their behavior for Bernoulli shifts, including estimates for the canonical generating partitions:
\begin{theorem*}[Theorem \ref{theorem: Canonical Bernoulli calculations}: Canonical Bernoulli calculations]

The canonical partition $\mathcal C$ is admissible (see Definition \ref{def: Return-Time Admissibility}). For every
$n\ge1$, we have 
\[n-\frac{\rho}{(1-\rho)^2}
 \le A_{\mu_{\mathbf p}}(\mathcal C_n)\le n \qquad \text{and} \qquad  V_{\mu_{\mathbf p}}(\mathcal C_n)
 \le \frac{\rho(1+\rho)}{(1-\rho)^3}.\]
Consequently,
\[L_{\mu_{\mathbf p}}(\sigma,\mathcal C)
=\frac1{H(\mathbf p)} \qquad \text{and} \qquad L_{V,\mu_{\mathbf p}}(\sigma,\mathcal C)=0,\]
and hence
\[\Lambda_{\mu_{\mathbf p}}(\sigma)\ge\frac1{H(\mathbf p)} \qquad \text{and} \qquad \Lambda_{V,\mu_{\mathbf p}}(\sigma)\ge0.\]
Here, $\textbf{p}$ is the probability distribution vector associate with the Bernoulli shift.
\end{theorem*}

\noindent For finite measure preserving periodic systems, we obtain explicit formulas as in the topological setting:

\begin{theorem*}[Theorem \ref{Theorem: Exact global values for finite permutations}: Exact Global Values for Finite Permutations]
For a finite full-support probability-preserving permutation, we have
\[ \Lambda_\mu(f)= \frac{\displaystyle \sum_{i=1}^r a_i m_i}
{\displaystyle 1-\sum_{i=1}^r a_i\log\left( \frac{a_i}{m_i} \right)} \qquad \text{and} \qquad 
\Lambda_{V,\mu}(f)=
\frac{\displaystyle\sum_{i=1}^r a_i
\left(m_i-\sum_{j=1}^r a_jm_j\right)^2}
{\displaystyle \left(1-\sum_{i=1}^r a_i\log\left(\frac{a_i}{m_i} \right)\right)^2}.\]
\end{theorem*}

\noindent Finally, we consider finite periodic approximations
of the uniform Bernoulli shift on $q$ symbols. For the corresponding finite systems $(X_N,f_N,\mu_N)$, we prove:

\begin{theorem*}[Theorem \ref{Theorem:Periodic Orbit Approximation}: Periodic Orbit Approximation]
For the finite systems defined in \eqref{Finite Systems}, we have
\begin{align*}
    \lim_{N\to\infty}\Lambda_{\mu_N}(f_N) &= \lim_{N\to\infty} \frac{\widetilde d_N}{1+\log s_N} = \frac{1}{\log q}, \quad \text{and}\\
    \\
    \lim_{N\to\infty}\Lambda_{V,\mu_N}(f_N) &= \lim_{N\to\infty}
\frac{\displaystyle \frac{1}{s_N}\sum_{d=1}^{N} P_q(d)(d-\widetilde d_N)^2 }{ (1+\log s_N)^2} = 0.
\end{align*}
\end{theorem*}

\noindent These results illustrate how normalized local return-time statistics can exhibit nontrivial asymptotic behavior even along sequences of finite systems that have zero entropy.

\subsection{Organization of The Paper}

The paper is organized as follows:

\noindent Section \ref{Section: Periodic Penetration} introduces
periodic penetration times, finite generating covers, and the asymptotic mean invariant $\Lambda_A$. We establish conjugacy invariance and study finite periodic systems, full shifts, shifts of finite type, and a zero entropy example.

\noindent Section \ref{Topological Invertible Setting} extends the construction to invertible systems using symmetric generating covers and two-sided dynamical refinements.

\noindent Section \ref{sec:variance} introduces the topological
variance invariant $\Lambda_V$ and investigates its properties in finite and symbolic dynamical systems, including the invertible setting.

\noindent Section \ref{section: Measure-Theoretic Return-Time Averages} develops the measure theoretic framework based on positive measure overlap times. We define the corresponding mean and variance invariants, establish their isomorphism invariance, and study finite systems, Bernoulli shifts, and periodic approximations.

\noindent Finally, Section \ref{Section: Further Questions} presents several open questions concerning admissibility, global invariant values, and relationships between the topological and measure theoretic constructions.

\section{Periodic Penetration and the Invariant $\Lambda_A$}\label{Section: Periodic Penetration}

Let $X$ be a non-empty compact Hausdorff space, and let $\displaystyle f:X\longrightarrow X$ be continuous. Throughout this section, we assume that the periodic points of $f$ are dense in $X$, that is, $\displaystyle \overline{\operatorname{Per}(f)}=X.$ If the set of all periodic points is not dense in $X$, then we restrict $f$ to $\overline{\operatorname{Per}(f)}$. In both cases, we assume $\overline{\operatorname{Per}(f)}\neq \emptyset$.

\noindent For a periodic point $x\in X$, let $\operatorname{per}_f(x)$ denote its
least period.

\begin{definition}[Periodic penetration time]
Let $U\subseteq X$ be a nonempty open set. Define
\[
\eta_f(U)
:=
\min\left\{
\operatorname{per}_f(x):
x\in U\cap\operatorname{Per}(f)
\right\}.
\]
We call $\eta_f(U)$ the \emph{periodic penetration time} of $U$.
\end{definition}
Since $\operatorname{Per}(f)$ is dense in $X$, every nonempty open set contains a periodic point. Hence, $0<\eta_f(U)<\infty$. Notice that if $V\subseteq U$, then $\displaystyle
\eta_f(V)\geq \eta_f(U)$.

\begin{proposition} Let $X$ and $Y$ be non-empty compact Hausdorff spaces. Let $\displaystyle f:X\longrightarrow X$ and $\displaystyle g:Y\longrightarrow Y$ be continuous. 
Suppose that
\[
h:(X,f)\longrightarrow (Y,g)
\]
is a topological conjugacy. Thus, $\displaystyle h\circ f (x) = g\circ h(x)$ for all $x\in X$. Then, for every nonempty open set
$U\subseteq X$,
\[
\eta_g(h(U))=\eta_f(U).
\]
\end{proposition}

\begin{proof}
A topological conjugacy preserves least periods. Indeed,
\[
f^n(p)=p
\quad\text{if and only if} \quad
g^n(h(p))=h(p).
\]
Thus, $p$ and $h(p)$ have the same least period. Since $h$ gives a
bijection from $U\cap\operatorname{Per}(f)$ onto
$h(U)\cap\operatorname{Per}(g)$, the result follows immediately.
\end{proof}

Throughout this paper, we assume $X$ to be a non-empty compact Hausdorff space.
\begin{definition}
Let
\[
\operatorname{Cov}_{\mathrm{fin}}(X)
:=
\left\{
\mathcal U :
\mathcal U \text{ is a finite open cover of }X
\right\}
\] be the collection of all finite open covers of $X$.
For $\displaystyle
\mathcal U\in \operatorname{Cov}_{\mathrm{fin}}(X)$,
define
\[
A_f(\mathcal U)
:=
\frac{1}{|\mathcal U|}
\sum_{U\in\mathcal U}\eta_f(U)
\]
and
\[
M_f(\mathcal U)
:=
\max_{U\in\mathcal U}\eta_f(U).
\]
\end{definition}
\noindent Thus, $A_f(\mathcal U)$ measures the average periodic penetration time
of the members of $\mathcal U$, whereas $M_f(\mathcal U)$ measures the
largest periodic penetration time among them. In particular,
\[
A_f(\mathcal U)\leq M_f(\mathcal U).
\]

\noindent For a finite open cover $\mathcal U$, let
\[
N(\mathcal U)
:=
\min\left\{
|\mathcal V|:
\mathcal V\subseteq\mathcal U
\text{ and }
\bigcup_{V\in\mathcal V}V=X
\right\}
\]
denote the minimum cardinality of a subcover of $\mathcal U$.

\begin{definition}[Dynamical join]\label{Dynamical Join}
For $\displaystyle
\mathcal U\in \operatorname{Cov}_{\mathrm{fin}}(X)$, and $n\geq1$, define
\[
\mathcal U_f^{(n)}
:=
\bigvee_{j=0}^{n-1}f^{-j}\mathcal U.
\]
Thus, $\mathcal U_f^{(n)}$ consists of the nonempty sets of the form
\[ U_{i_0}\cap f^{-1}(U_{i_1}) \cap\cdots\cap f^{-(n-1)}(U_{i_{n-1}}), \qquad \text{where} \qquad U_{i_j}\in\mathcal U. \]
\end{definition}

\begin{definition}[Finite topological generating cover]
Let $\mathcal U$ be a finite open cover of $X$. We say that
$\mathcal U$ is a \emph{finite topological generating cover} for
$(X,f)$ if, for every
$\mathcal V\in\operatorname{Cov}_{\mathrm{fin}}(X)$, there exists
$n\geq 1$ such that $\displaystyle \mathcal U_f^{(n)}\prec\mathcal V$.

\noindent Here $\mathcal U_f^{(n)}\prec\mathcal V$ means that
$\mathcal U_f^{(n)}$ refines $\mathcal V$; that is, for every
$U\in\mathcal U_f^{(n)}$, there exists $V\in\mathcal V$ such that $U\subseteq V.$ We denote the collection of all finite topological generating covers
for $(X,f)$ by
\[
\operatorname{Gen}(X,f).
\]
\end{definition}

\begin{definition}[Generating cardinality]
Suppose that $\operatorname{Gen}(X,f)\neq\varnothing$. Define 

\noindent $\displaystyle r(f): \operatorname{Gen}(X,f) \rightarrow \mathbb{N}_0$ by
\[ r(f) := \min\left\{ |\mathcal U|: \mathcal{U} \in\operatorname{Gen}(X,f) \right\}. \]
\end{definition}

\begin{proposition}\label{generating cover proposition}
Let $(X,f)$ and $(Y,g)$ be topologically conjugate dynamical systems, and let $\displaystyle h:X\longrightarrow Y $ be a topological conjugacy. Then,
\[ \mathcal U\in\operatorname{Gen}(X,f) \quad \text{if and only if} \quad
h(\mathcal U)\in\operatorname{Gen}(Y,g), \]
where
\[ h(\mathcal U):=\{h(U):U\in\mathcal U\}. \]
Consequently,
\[ r(f)=r(g). \]
\end{proposition}

\begin{proof}
Since $h$ is a homeomorphism, $h(\mathcal U)$ is a finite open cover of $Y$ whenever $\mathcal U$ is a finite open cover of $X$.

\noindent $(\Rightarrow)$
Suppose $\displaystyle \mathcal U\in\operatorname{Gen}(X,f)$.
Let $\mathcal V\in\operatorname{Cov}_{\mathrm{fin}}(Y)$. Since $h$ is a homeomorphism, we get
\[
h^{-1}(\mathcal V)
:=
\{h^{-1}(V):V\in\mathcal V\}
\]
is a finite open cover of $X$. Thus, by the definition of $\operatorname{Gen}(X,f)$, there exists $n\geq 1$ such that
$\displaystyle \mathcal U_f^{(n)}\prec h^{-1}(\mathcal V)$.
\noindent Since $h:X \longrightarrow Y$ is a bijection, and $h\circ f=g\circ h$, we have $\displaystyle h\bigl(\mathcal U_f^{(n)}\bigr) = h(\mathcal U)_g^{(n)}$.
Therefore, $\displaystyle h(\mathcal U)_g^{(n)}\prec\mathcal V$. Thus, $\displaystyle h(\mathcal U)\in\operatorname{Gen}(Y,g)$.

\noindent $(\Leftarrow)$ By applying the same argument to the inverse conjugacy $\displaystyle h^{-1}:(Y,g)\longrightarrow(X,f)$
we get the converse implication. 

\noindent Hence,
\[
\mathcal U\in\operatorname{Gen}(X,f)
\quad \text{if and only if} \quad
h(\mathcal U)\in\operatorname{Gen}(Y,g).
\]

\noindent Moreover, since $h$ is injective, we have $\displaystyle |h(\mathcal U)|=|\mathcal U|$. Thus, $h$ induces a cardinality-preserving bijection between
$\operatorname{Gen}(X,f)$ and $\operatorname{Gen}(Y,g)$. This implies, $\displaystyle r(f)=r(g)$.
\end{proof}

\begin{definition}[$\eta$-admissible generating cover]\label{Definition: eta admissible generating covers}
Let $\mathcal U\in\operatorname{Gen}(X,f)$. We say that
$\mathcal U$ is \emph{$\eta$-admissible} if
\[ \limsup_{n\to\infty} \frac{M_f\bigl(\mathcal U_f^{(n)}\bigr)}{n} <\infty.\]

\noindent The admissibility condition requires periodic orbit penetration to keep pace, at most linearly, with the dynamical refinement generated by $\mathcal U$. More precisely, while the iterated joins $\mathcal U_f^{(n)}$ resolve increasingly fine orbit information, $\eta$-admissibility requires the largest least period needed to penetrate every member of $\mathcal U_f^{(n)}$ to grow at most linearly in $n$. Thus, admissibility excludes generating covers for which periodic orbit penetration becomes disproportionately slow relative to the length of the observed orbit segment.

\noindent We denote the collection of all $\eta$-admissible finite topological generating covers by
\[ \mathscr A_f := \left\{\mathcal {U}\in\operatorname{Gen}(X,f):
\limsup_{n\to\infty} \frac{M_f\bigl(\mathcal U_f^{(n)}\bigr)}{n}
<\infty \right\}.\]
\end{definition}

\begin{proposition}\label{conjugacy preserves non-empty}
Let $(X,f)$ and $(Y,g)$ be topologically conjugate dynamical systems, and let $\displaystyle h:X\longrightarrow Y$
be a topological conjugacy. Then,
\[\mathcal U\in\mathscr A_f \quad \text{if and only if} \quad h(\mathcal U)\in\mathscr A_g.\]
In particular,
\[\mathscr A_f\neq\varnothing \quad \text{if and only if} \quad
\mathscr A_g\neq\varnothing.\]
\end{proposition}

\begin{proof}
Suppose that $\mathcal U\in\mathscr A_f$. Then
$\mathcal U\in\operatorname{Gen}(X,f)$. By Proposition
\ref{generating cover proposition},
\[
h(\mathcal U)\in\operatorname{Gen}(Y,g).
\]
Moreover, for every $n\geq1$,
\[
h\bigl(\mathcal U_f^{(n)}\bigr)
=
h(\mathcal U)_g^{(n)}.
\]
Since topological conjugacy preserves least periods of periodic points,
\[
M_g\bigl(h(\mathcal U)_g^{(n)}\bigr)
=
M_f\bigl(\mathcal U_f^{(n)}\bigr).
\]
Therefore,
\[
\limsup_{n\to\infty}
\frac{M_g\bigl(h(\mathcal U)_g^{(n)}\bigr)}{n}
=
\limsup_{n\to\infty}
\frac{M_f\bigl(\mathcal U_f^{(n)}\bigr)}{n}
<\infty.
\]
Hence,
\[
h(\mathcal U)\in\mathscr A_g.
\]

\noindent Applying the same argument to the inverse conjugacy
$h^{-1}:(Y,g)\longrightarrow(X,f)$ gives the converse implication.
Therefore,
\[
\mathcal U\in\mathscr A_f
\quad \text{if and only if} \quad
h(\mathcal U)\in\mathscr A_g.
\]
In particular,
\[
\mathscr A_f\neq\varnothing
\quad \text{if and only if} \quad
\mathscr A_g\neq\varnothing.
\]
\end{proof}

\begin{definition}
For $\mathcal U\in\mathscr A_f$, define
\[
L_A(f,\mathcal U)
:=
\limsup_{n\to\infty}
\frac{
A_f\bigl(\mathcal U_f^{(n)}\bigr)
}{
1+\log N\bigl(\mathcal U_f^{(n)}\bigr)
}.
\]
\end{definition}

\noindent The quantity $L_A(f,\mathcal U)$ is associated with the pair $(f,\mathcal U)$. Thus, it may depend on the choice of the generating cover. We therefore define a system-level quantity by taking the supremum over all $\eta$-admissible generating covers.

\begin{definition}\label{Definition of Lambda_A}
Suppose $\mathscr A_f\neq\varnothing$. Define the system-level supremal asymptotic mean periodic penetration rate
\[\Lambda_A(f) := \sup_{\mathcal U\in\mathscr A_f} L_A(f,\mathcal U).\]
\end{definition}

\noindent For a fixed admissible generating cover $\mathcal U$, the quantity $L_A(f,\mathcal U)$ measures the asymptotic average periodic penetration of the dynamically refined covers, normalized by the logarithm of their minimum subcover cardinalities. The invariant $\Lambda_A(f)$ takes the supremum of this quantity over all $\eta$-admissible generating covers. Thus, $\Lambda_A(f)$ records the largest normalized average periodic penetration rate detectable through an admissible generating cover of the system.

\begin{theorem}[Conjugacy Invariance]\label{Theorem: Lambda is Conjugacy Invariant}
Suppose that $\displaystyle h:(X,f)\longrightarrow(Y,g)$ is a topological conjugacy, and suppose that the preceding hypotheses hold for $(X,f)$ and $(Y,g)$. Then,
\[\Lambda_A(f)=\Lambda_A(g).\]
\end{theorem}

\begin{proof}
Let $\mathcal U\in\mathscr A_f$. Since $h$ is a topological
conjugacy, for every $n\geq1$, $\displaystyle h\bigl(\mathcal U_f^{(n)}\bigr) = h(\mathcal U)_g^{(n)}$.
Since conjugacy preserves periodic penetration times, we have
$\displaystyle A_f\bigl(\mathcal U_f^{(n)}\bigr) = A_g\bigl(h(\mathcal U)_g^{(n)}\bigr)$. Moreover, since $h$ induces a bijection between the members of $\mathcal U_f^{(n)}$ and those of $h(\mathcal U)_g^{(n)}$ and preserves inclusions, the corresponding subcovers have the same cardinality. Hence,
$\displaystyle N\bigl(\mathcal U_f^{(n)}\bigr) = N\bigl(h(\mathcal U)_g^{(n)}\bigr)$. Therefore,
\[ L_A(f,\mathcal U) = L_A(g,h(\mathcal U)).\]

\noindent By Proposition \ref{conjugacy preserves non-empty}, the map $\displaystyle \mathcal U\longmapsto h(\mathcal U)$ is a bijection from $\mathscr A_f$ onto $\mathscr A_g$. Taking suprema over the corresponding admissible generating covers gives
\[\Lambda_A(f)=\Lambda_A(g).\]
\end{proof}

\subsection{Topological Entropy and Generating Covers}

We briefly recall the open-cover definition of topological entropy introduced by Adler, Konheim, and McAndrew in 1965 \cite{topological_entropy_adler1965}. For alternative definitions and further details, one may refer to standard references such as \cite{entropy_bowen1971}, \cite{katok_and_Hasselblatt_book_1995introduction}, and \cite{walters_book_1982introduction}.

\noindent Let $\mathcal U\in\operatorname{Cov}_{\mathrm{fin}}(X)$. The entropy of $f$ relative to $\mathcal U$ is defined by
\[ h(f,\mathcal U) := \lim_{n\to\infty} \frac{1}{n} \log N\!\left(\mathcal U_f^{(n)}\right), \]
where
\[ \mathcal U_f^{(n)} = \bigvee_{j=0}^{n-1}f^{-j}\mathcal U. \]
The limit exists by the subadditivity of $\displaystyle \left(\log N(\mathcal U_f^{(n)})\right)_{n=1}^\infty$ due to Fekete's Lemma. The topological entropy of $f$ is then
defined by \[ h_{\mathrm{top}}(f) := \sup_{\mathcal V\in\operatorname{Cov}_{\mathrm{fin}}(X)} h(f,\mathcal V). \]

\begin{lemma}\label{Entropy Lemma}
If $\mathcal U\in\operatorname{Gen}(X,f)$, then
\[ h(f,\mathcal U)=h_{\mathrm{top}}(f). \]
\end{lemma}

\begin{proof}
Let $\mathcal V\in\operatorname{Cov}_{\mathrm{fin}}(X)$. Since
$\mathcal U$ is generating, there exists $m\geq 1$ such that
$\displaystyle \mathcal U_f^{(m)}\prec\mathcal V$.
Consequently, for every $n\geq 1$, we have $\displaystyle \mathcal U_f^{(m+n-1)} \prec \mathcal V_f^{(n)}$.
Hence, $\displaystyle N\!\left(\mathcal V_f^{(n)}\right) \leq N\!\left(\mathcal{U}_f^{(m+n-1)}\right)$.
Therefore, $\displaystyle h(f,\mathcal V) \leq h(f,\mathcal U)$.
Since $\mathcal V$ was arbitrary, we get $\displaystyle  h_{\mathrm{top}}(f)\leq h(f,\mathcal U)$.
The reverse inequality follows immediately from the definition of
topological entropy. Thus,
\[ h(f,\mathcal U)=h_{\mathrm{top}}(f). \]
\end{proof}

\subsection{Finite Periodic Orbits}

We begin with a finite system for which $\Lambda_A$ can be computed explicitly.

\begin{proposition}\label{calculating Lambda for finite systmes}
Let $X$ be a finite dynamical system consisting of $r$ distinct
periodic orbits with least periods $m_1,\ldots,m_r$, and let
$f:X\to X$ be the corresponding map. Then,
\[\Lambda_A(f) = \frac{\displaystyle\sum_{i=1}^{r}m_i^2}
{\displaystyle \left(\sum_{i=1}^{r}m_i\right) \left(1+\log\left(\sum_{i=1}^{r}m_i\right)\right)}.\]
\end{proposition}

\begin{proof}
Since $X$ is finite and Hausdorff, it has the discrete topology.
Let $\mathcal U$ be any finite generating cover of $X$. By the generating property, there exists $N\geq1$ such that
\[
\mathcal U_f^{(N)}\prec
\bigl\{\{x\}:x\in X\bigr\}.
\]
Since iterated joins become finer as $n$ increases, for every
$n\geq N$, the distinct nonempty members of $\mathcal U_f^{(n)}$
are precisely the singletons of $X$. For every $x\in X$, we have
\[\eta_f(\{x\})=\operatorname{per}_f(x).\]
Consequently, for all $n\geq N$, $\displaystyle M_f(\mathcal U_f^{(n)}) = \max_{1\leq i\leq r}m_i$. Thus, every finite generating cover is $\eta$-admissible. Moreover,
\[N(\mathcal U_f^{(n)}) = |X| = \sum_{i=1}^{r}m_i. \]
Since the $i$-th periodic orbit contains exactly $m_i$ points,
each having least period $m_i$, we obtain
\[A_f(\mathcal U_f^{(n)}) = \frac{1}{|X|} \sum_{x\in X}\operatorname{per}_f(x) = \frac{\displaystyle \sum_{i=1}^{r}m_i^2} {\displaystyle \sum_{i=1}^{r}m_i}.\]

\noindent Therefore,
\[ L_A(f,\mathcal U) = \frac{\displaystyle\sum_{i=1}^{r}m_i^2}
{\displaystyle \left(\sum_{i=1}^{r}m_i\right) \left(1+\log\left(\sum_{i=1}^{r}m_i\right)\right)}.\]

\noindent Since this value is independent of the generating cover
$\mathcal U$, taking the supremum over all $\eta$-admissible
generating covers gives the desired formula.
\end{proof}

\begin{remark}
The preceding proposition shows that $\Lambda_A$ contains
information not detected by topological entropy. Indeed, every
finite dynamical system has zero topological entropy, whereas
\[\Lambda_A(f) = \frac{\displaystyle\sum_{i=1}^{r}m_i^2}
{\displaystyle\left(\sum_{i=1}^{r}m_i\right)\left(1+\log\left(\sum_{i=1}^{r}m_i\right)\right)}\]
depends on the lengths of its periodic orbits. In particular, for a single periodic orbit of length $m$,
\[\Lambda_A(f)=\frac{m}{1+\log m}.\]
Thus, even among zero-entropy systems, $\Lambda_A$ can distinguish
different periodic orbit structures.
\end{remark}

\begin{proposition}
Let $(\Sigma_q^+,\sigma_q)$ be the one-sided full shift on
$q\geq2$ symbols. For each $N\geq1$, define
\[X_N = \bigcup_{k=1}^{N}\operatorname{Fix}(\sigma_q^k),
\qquad f_N=\sigma_q|_{X_N}.\]
Then,
\[\lim_{N\to\infty}\Lambda_A(f_N) = \frac{1}{\log q}. \]
\end{proposition}

\begin{proof}
For each $d\geq1$, let $P_q(d)$ denote the number of points
of least period exactly $d$ under $\sigma_q$.  Since every point of $X_N$ has least period at most $N$, we have
\[|X_N|=\sum_{d=1}^{N}P_q(d).\]
By the preceding proposition,
\[\Lambda_A(f_N) = \frac{\displaystyle\sum_{d=1}^{N}dP_q(d)}
{\displaystyle \left(\sum_{d=1}^{N}P_q(d)\right) \left(1+\log\left(\sum_{d=1}^{N}P_q(d)\right)\right)}.\]

\noindent Notice that, for the full shift, the number of points fixed by $\sigma_q^d$ is $q^d$. Hence,
\[q^d=\sum_{k\mid d}P_q(k).\]
By Möbius inversion,
\[P_q(d) = \sum_{k\mid d}\mu(k)q^{d/k},\]
where $\mu$ denotes the Möbius function. Since every point of least period $d$ is fixed by $\sigma_q^d$, we have
\[P_q(d)\leq |\operatorname{Fix}(\sigma_q^d)|=q^d.\]

\noindent In particular,
\[P_q(d)=q^d+O(dq^{d/2}).\]
Consequently,
\[\sum_{d=1}^{N}P_q(d) = \sum_{d=1}^{N}q^d+O(N^2q^{N/2})
= \frac{q^{N+1}-q}{q-1}+O(N^2q^{N/2}). \]
It follows that
\[\sum_{d=1}^{N}P_q(d) \sim \frac{q^{N+1}}{q-1},\]
and therefore
\[\log\left(\sum_{d=1}^{N}P_q(d)\right) = N\log q+O(1).\]

\noindent Furthermore,
\[0\leq N- \frac{\displaystyle \sum_{d=1}^{N}dP_q(d)}{\displaystyle \sum_{d=1}^{N}P_q(d)} = \frac{\displaystyle \sum_{d=1}^{N}(N-d)P_q(d)} {\displaystyle \sum_{d=1}^{N}P_q(d)}.\]
Since $P_q(d)\leq q^d$, the numerator satisfies
\[\sum_{d=1}^{N}(N-d)P_q(d) \leq \sum_{d=1}^{N}(N-d)q^d = q^N\sum_{j=0}^{N-1}\frac{j}{q^j} = O(q^N).\]
Because $\displaystyle \sum_{d=1}^{N}P_q(d)\sim\frac{q^{N+1}}{q-1}$, we have 

\[\frac{\displaystyle \sum_{d=1}^{N}dP_q(d)}{\displaystyle \sum_{d=1}^{N}P_q(d)} = N+O(1) \qquad \text{for sufficiently large enough } N.\]

\noindent Therefore,
\[\lim_{N \rightarrow \infty}\Lambda_A(f_N) = \lim_{N \rightarrow \infty}\frac{N+O(1)}{1+N\log q+O(1)} = \frac{1}{\log q}.\]
\end{proof}

\begin{remark}
The finite invariant sets $X_N$ form an increasing sequence,
and $\displaystyle \overline{\bigcup_{N\geq 1}X_N}=\Sigma_q^+$,
since periodic points are dense in the full shift. The preceding proposition shows that the periodic penetration invariants of these finite subsystems converge to $\displaystyle \frac{1}{\log q}$. As we shall see in the next subsection, this is
also the value obtained from the standard cylinder generating
partition of the full shift. However, this convergence does not by itself imply that
\[\Lambda_A(\sigma_q)=\frac{1}{\log q},\]
as the supremum defining $\Lambda_A(\sigma_q)$ ranges over
all $\eta$-admissible generating covers.
\end{remark}

\subsection{The One-Sided Full Shift}

We now compute the invariant $\Lambda_A$ for the one-sided full shift. For more details about symbolic spaces, one may refer to classical texts such as \cite{kitchens1998symbolic}, \cite{lind_marcus_1995}.

\noindent Let
\[
\Sigma_q^+
=
\{0,1,\ldots,q-1\}^{\mathbb N_0},
\qquad q\geq 2,
\]
be the one-sided sequence space on $q$-symbols, and let
$\displaystyle \sigma_q:\Sigma_q^+\longrightarrow\Sigma_q^+$
be the left shift defined by
\[
(\sigma_q(x))_j=x_{j+1}.
\]

\noindent For a word
\[
w=w_0w_1\cdots w_{n-1} \in\{0,1,\ldots,q-1\}^n, \]
of length $n$, the corresponding cylinder set is
\[ [w]_k = \left\{ x=(x_j)_{j\geq0}\in\Sigma_q^+: x_{k+i}=w_i,\ 0\leq i \leq n-1 \right\}. \]
For simplicity, we let $[w] \equiv [w]_0$.

\noindent We denote the collection of all $n$-cylinders by
\[
\mathcal C_n
=
\left\{
[w]:
w\in\{0,1,\ldots,q-1\}^n
\right\}.
\]

\noindent Under the metric on $\Sigma^+_n$ defined by $\displaystyle d(x,y):=\frac{1}{2^M}$ where $\displaystyle M:=\max \{k: x_i=y_i\ \text{for}\ i<k\}$, the collection $\mathcal C_n$ is a clopen partition of $\Sigma_q^+$. Furthermore, we have $\displaystyle |\mathcal C_n|=q^n$.

\begin{proposition}
For every $n\geq1$, we have $\displaystyle (\mathcal C_1)_{\sigma_q}^{(n)} = \mathcal C_n$.  Consequently,
$\displaystyle N(\mathcal C_n)=q^n$.
\end{proposition}

\begin{proof}
By definition,
\[
(\mathcal C_1)_{\sigma_q}^{(n)}
=
\bigvee_{j=0}^{n-1}\sigma_q^{-j}\mathcal C_1.
\]
An element of this join is of the form $\displaystyle [i_0]
\cap\sigma_q^{-1}[i_1]
\cap\cdots\cap
\sigma_q^{-(n-1)}[i_{n-1}]$. This is precisely the cylinder
$\displaystyle [i_0i_1\cdots i_{n-1}$.
Hence,
\[
(\mathcal C_1)_{\sigma_q}^{(n)}=\mathcal C_n.
\]

\noindent Since the elements of $\mathcal C_n$ are pairwise disjoint and form a
partition of $\Sigma_q^+$, every subcover must contain every member of
$\mathcal C_n$. Therefore,
\[
N(\mathcal C_n)=|\mathcal C_n|=q^n.
\]
\end{proof}

\noindent Next, we determine the periodic penetration time of a cylinder. For a word $\displaystyle w=w_0w_1\cdots w_{n-1}$,
define its least word-period by
\[
p(w)
:=
\min\left\{
1\leq p\leq n:
w_i=w_{i+p}
\text{ whenever }0\leq i<i+p<n
\right\}.
\]

\begin{proposition}
For each $n\in \mathbb{N}$, and for every word $w$ of length $n$, we have
\[
\eta_{\sigma_q}([w])=p(w).
\]
\end{proposition}

\begin{proof}
Let $p\equiv p(w)$. By the definition of $p(w)$, the word $w$ is compatible
with repetition with period $p$. Thus, the periodic sequence
$\displaystyle x=(w_0w_1\cdots w_{p-1})^\infty$ belongs to $[w]$. Moreover, the least period of $x$ is $p$; otherwise, 
$w$ would have a smaller word-period. Hence,
\begin{align}\label{eta([w])<p(w)}
    \eta_{\sigma_q}([w])\leq p(w).
\end{align}

\noindent Conversely, suppose that $x\in[w]$ is periodic with least period $r$. Then, $\displaystyle x_i=x_{i+r}$ whenever both coordinates occur among the first $n$ coordinates.
Since the first $n$ coordinates of $x$ are given by $w$, we get $w_i=w_{i+r}$ whenever $0\leq i<i+r<n$. Therefore, $r$ is a period of $w$, and hence $p(w)\leq r$.
Taking the minimum over all periodic $x\in[w]$ gives
\begin{align}\label{p(w)<eta([w])}
    p(w)\leq\eta_{\sigma_q}([w]).
\end{align}
Thus, by (\ref{eta([w])<p(w)}) and (\ref{p(w)<eta([w])}), we get
\[
\eta_{\sigma_q}([w])=p(w).
\]
\end{proof}

\begin{corollary}
For every $n\geq1$, $\displaystyle M_{\sigma_q}(\mathcal C_n)=n$.
\end{corollary}

\begin{proof}
Since every word of length $n$ has a word period of at most $n$, we have $\displaystyle \eta_{\sigma_q}([w])\leq n$ for every $[w]\in\mathcal C_n$. Hence, $\displaystyle M_{\sigma_q}(\mathcal C_n)\leq n$. On the other hand, there exist words of length $n$ whose least
word-period is $n$. For example,
\[
0^{n-1}1\equiv \underbrace{0\cdots 0}_{n\text{-times}}1
\]
has least word-period $n$. That is; $\displaystyle p\left(0^{n-1}1\right)=n$. Therefore,
\[
M_{\sigma_q}(\mathcal C_n)=n.
\]
\end{proof}

\noindent We now estimate the average periodic penetration time of the
$n$-cylinder partition.

\begin{proposition}\label{A(C_n)/n rightarrow 1}
For the one-sided full $q$-shift $\Sigma_q^+$,
\[
0\le
n-A_{\sigma_q}(\mathcal C_n)
\le
\frac{q}{(q-1)^2}
\qquad
\text{for every } n\in\mathbb N.
\]
Consequently,
\[
\lim_{n\to\infty}
\frac{A_{\sigma_q}(\mathcal C_n)}{n}=1.
\]
\end{proposition}

\begin{proof}
Let $n\in \mathbb{N}$. Because $\displaystyle \eta_{\sigma_q}([w])=p(w)$, we have
\[ A_{\sigma_q}(\mathcal C_n) = \frac{1}{q^n} \sum_{|w|=n}p(w).\]

\noindent Recall: $\displaystyle \sum_{|w|=n}$ means the sum over all words $w$ of length $n$. Because $\displaystyle \mathcal{C}_n = \{[w]: |w|=n\}$, one may also use $\displaystyle A_{\sigma_q}(\mathcal {C}_n) = \frac{1}{|\mathcal{C}_n|}
\sum_{[w]\in \mathcal{C}_n}p(w)$.\\

\noindent Set $ D_n(w):=n-p(w)$. Then,
\[ n-A_{\sigma_q}(\mathcal C_n) = \frac{1}{q^n}\sum_{|w|=n}D_n(w). \]

\noindent For $1\leq k\leq n-1$, if $ D_n(w)\geq k$, then $
p(w)\leq n-k$. Notice that for each $1\leq r\leq n-k$, a word of length $n$ having period $r$ is completely determined by its first $r$ symbols. Hence, there are at most $q^r$ such words. Therefore,
\[ \left| \{w: |w|=n, D_n(w)\geq k\} \right| \leq \sum_{r=1}^{n-k}q^r
\leq \frac{q^{n-k+1}}{q-1}. \]
Dividing by $q^n$ gives
\[ \frac{ |\{w: |w|=n, D_n(w)\geq k\}| }{q^n} \leq \frac{q^{1-k}}{q-1}. \]

\noindent Since $D_n(w)\in\{0,1,\ldots,n-1\}$, we have
\[ D_n(w) = \sum_{k=1}^{n-1} \mathbf{1}_{\{D_n(w)\ge k\}}. \]
Therefore, summing over all words $w$ of length $n$ and interchanging the finite sums gives
\begin{align*}
\frac{1}{q^n}\sum_{|w|=n}D_n(w)
&= \frac{1}{q^n} \sum_{|w|=n} \sum_{k=1}^{n-1} \mathbf{1}_{\{D_n(w)\ge k\}}\\
&= \sum_{k=1}^{n-1} \frac{ \left|\{w:|w|=n,\ D_n(w)\ge k\}\right| }{q^n}.
\end{align*}

\noindent Hence,
\[ 0\leq n-A_{\sigma_q}(\mathcal C_n) \leq \sum_{k=1}^{\infty} \frac{q^{1-k}}{q-1} = \frac{q}{(q-1)^2}\quad \text{for every}\quad n\in \mathbf{N}.\]
Consequently,
\[ \lim_{n\rightarrow \infty}\frac{A_{\sigma_q}(\mathcal C_n)}{n} = 1.\]
\end{proof}

\noindent We next determine the generating cardinality of the one-sided full shift.

\begin{proposition} \label{r(sigma_q)=q}
For the one-sided full $q$-shift, we have
\[
r(\sigma_q)=q.
\]
\end{proposition}

\begin{proof}
First, recall that $\displaystyle r(\sigma_q)= \left\{|\mathcal{U}|: \mathcal{U}\in \operatorname{Gen} \left(\Sigma_q^+,\sigma_q \right) \right\}$. The one-symbol cylinder partition
\[ \mathcal C_1=\{[0],[1],\ldots,[q-1]\} \]
is a finite topological generating cover because $\displaystyle (\mathcal C_1)_{\sigma_q}^{(n)}=\mathcal C_n$, the cylinder partitions become arbitrarily fine. Hence,
\begin{align}\label{r(sigma_q)<= q}
    r(\sigma_q)\leq q.
\end{align}

\noindent For the reverse inequality, let
$\mathcal U\in\operatorname{Gen}(\Sigma_q^+,\sigma_q)$ and suppose
$\displaystyle |\mathcal U|=m$.
By Lemma \ref{Entropy Lemma}, we get $\displaystyle h_{\mathrm{top}}(\sigma_q) = h(\sigma_q,\mathcal U)$. 
Moreover, we have $\displaystyle N\!\left(\mathcal {U}_{\sigma_q}^{(n)}\right) \leq \left|\mathcal{U}_{\sigma_q}^{(n)}\right| \leq m^n$. Therefore,
\[ h_{\mathrm{top}}(\sigma_q) = \lim_{n\to\infty} \frac{1}{n}
\log N\!\left(\mathcal U_{\sigma_q}^{(n)}\right) \leq \log m. \]
Since the one-sided full $q$-shift satisfies $\displaystyle h_{\mathrm{top}}(\sigma_q)=\log q$, we obtain $\displaystyle \log q\leq\log m$, and hence $m\geq q$. Since $\mathcal U$ was an arbitrary finite topological generating cover, it follows that
\begin{align}\label{r(sigma_q)>= q}
    r(\sigma_q)\geq q.
\end{align}
Combining (\ref{r(sigma_q)<= q}) and (\ref{r(sigma_q)>= q}), we get
\[r(\sigma_q)=q.\]
\end{proof}

\begin{theorem}[One-Sided Full Shifts]\label{L_A for One-Sided Full Shifts}
Let $\sigma_q:\Sigma_q^+\to\Sigma_q^+$ be the one-sided full shift
on $q\geq2$ symbols. Then the standard one-symbol cylinder partition $\mathcal C_1$ is $\eta$-admissible and
\[L_A(\sigma_q,\mathcal C_1) = \frac{1}{\log q}. \]
Consequently, 
\[ \Lambda_A(\sigma_q) \geq \frac{1}{\log q}. \]
\end{theorem}

\begin{proof}
Since $\displaystyle (\mathcal C_1)_{\sigma_q}^{(n)}=\mathcal C_n$
and $\displaystyle M_{\sigma_q}(\mathcal C_n)=n$, we have
\[ \limsup_{n\to\infty} \frac{M_{\sigma_q}(\mathcal C_n)}{n} =1<\infty. \]
Hence, $\displaystyle \mathcal C_1\in\mathscr A_{\sigma_q}$. Moreover,
\begin{align*}
L_A(\sigma_q,\mathcal C_1)
&= \limsup_{n\to\infty} \frac{A_{\sigma_q}(\mathcal C_n)}{1+\log N(\mathcal C_n)}\\
&= \limsup_{n\to\infty} \frac{A_{\sigma_q}(\mathcal C_n)/n} {1/n+\log q}\\
&= \frac{1}{\log q},
\end{align*}
where the last equality follows from Proposition
\ref{A(C_n)/n rightarrow 1}.

\noindent Therefore,
\[ \Lambda_A(\sigma_q) \geq L_A(\sigma_q,\mathcal C_1) = \frac{1}{\log q}. \]
\end{proof}

\begin{proposition}[Clopen generating covers of the full shift]\label{clopen cover proposition}
Let $\displaystyle \sigma_q:\Sigma_q^+\longrightarrow\Sigma_q^+$ 
be the one-sided full shift on $q\geq2$ symbols, and let $\mathcal U$ be a finite clopen generating cover of $\Sigma_q^+$. Then, $\mathcal U$ is $\eta$-admissible and
\[
L_A(\sigma_q,\mathcal U)
\leq
\frac{1}{\log q}.
\]
\end{proposition}

\begin{proof}
Since every member of $\mathcal U$ is clopen and $\mathcal U$ is
finite, there exists $K\geq1$ such that every $U\in\mathcal U$ is a
union of cylinders of length $K$. Let
\[V = U_0\cap\sigma_q^{-1}(U_1)\cap\cdots \cap\sigma_q^{-(n-1)}(U_{n-1}) \]
be a nonempty member of $\mathcal U_{\sigma_q}^{(n)}$. Let $x\in V$.
Since membership in each $U_j$ is determined by the first $K$
coordinates, every point agreeing with $x$ in the coordinates $0,1,\ldots,n+K-2$ also belongs to $V$. Let $p$ be the periodic point obtained by repeating the finite block $x_0x_1\cdots x_{n+K-2}$. Hence, $p\in V$, and $\operatorname{per}_{\sigma_q}(p)\leq n+K-1$.
Thus,
\[\eta_{\sigma_q}(V)\leq n+K-1.\]
Because $V$ was arbitrary, we get $\displaystyle M_{\sigma_q}\bigl(\mathcal U_{\sigma_q}^{(n)}\bigr) \leq n+K-1$. Therefore,
\[
\limsup_{n\to\infty}
\frac{
M_{\sigma_q}\bigl(\mathcal U_{\sigma_q}^{(n)}\bigr)
}{n}
\leq1,
\]
and hence $\mathcal U$ is $\eta$-admissible.

\noindent Furthermore,
\[ A_{\sigma_q}\bigl(\mathcal U_{\sigma_q}^{(n)}\bigr) \leq M_{\sigma_q}\bigl(\mathcal U_{\sigma_q}^{(n)}\bigr) \leq n+K-1.\]
Since $\mathcal U$ is generating, we have
\[\lim_{n\to\infty} \frac{1}{n} \log N\bigl(\mathcal U_{\sigma_q}^{(n)}\bigr) = h_{\mathrm{top}}(\sigma_q) = \log q. \]
Thus,
\begin{align*}
L_A(\sigma_q,\mathcal U)
&= \limsup_{n\to\infty}
\frac{ A_{\sigma_q}\bigl(\mathcal U_{\sigma_q}^{(n)}\bigr)
}{ 1+\log N\bigl(\mathcal U_{\sigma_q}^{(n)}\bigr) }\\
&\leq \limsup_{n\to\infty} \frac{n+K-1} {1+\log N\bigl(\mathcal U_{\sigma_q}^{(n)}\bigr)}\\
&= \frac{1}{\log q}.
\end{align*}
\end{proof}

\begin{corollary}
Among all finite clopen generating covers of the one-sided full shift,
the maximal value of $L_A$ is
\[
\frac{1}{\log q}.
\]
In particular, this value is attained by the standard one-symbol
cylinder partition $\mathcal C_1$.
\end{corollary}

\begin{proof}
By the preceding Proposition \ref{clopen cover proposition}, we have $\displaystyle L_A(\sigma_q,\mathcal U) \leq \frac{1}{\log q}$ for every finite clopen generating cover $\mathcal U$. On the other hand, the standard cylinder partition $\mathcal C_1$ satisfies
\begin{align}\label{L_A(sigma_q, C_1)}
    L_A(\sigma_q,\mathcal C_1) = \frac{1}{\log q}
\end{align}
\end{proof}

\begin{remark}
Since $\Lambda_A(\sigma_q)$ is defined using all $\eta$-admissible
finite open generating covers, the preceding result yields $\displaystyle  \Lambda_A(\sigma_q) \geq \frac{1}{\log q}$.
The question of whether equality holds when arbitrary open generating covers are allowed is not addressed here.
\end{remark}

\subsection{Mixing Shifts of Finite Type}\label{subsection: Mixing Shifts of Finite Type}
For completeness, we first present some basic ideas about shifts of
finite type (SFTs). For more details, one may refer to standard texts
such as \cite{kitchens1998symbolic}, \cite{weiss1973subshifts_sofic}. A \emph{subshift} is a closed,
shift-invariant subset of a full shift $\Sigma_q^+$ (or $\Sigma_q$).
A \emph{shift of finite type} is a subshift that can be described by
forbidding finitely many finite words. Thus, every shift of finite type is a subshift, although not every subshift is of finite type.

\noindent Let
\[ \mathcal {A}_q=\{0,1,\ldots,q-1\} \]
be a finite alphabet, and let
\[ B=(B_{ij})_{i,j\in\mathcal A} \]
be a $q\times q$ matrix with entries in $\{0,1\}$. We interpret
$B_{ij}=1$ to mean that the transition from the symbol $i$ to the
symbol $j$ is allowed, while $B_{ij}=0$ means that this transition is
forbidden. The associated one-sided shift space is
\[ \Sigma_B^+ = \left\{x=(x_0x_1x_2\ldots)\in\mathcal A^{\mathbb N_0}:
B_{x_i x_{i+1}}=1 \text{ for every }i\geq0 \right\}. \]

\noindent
A finite word $\displaystyle w=w_0w_1\cdots w_{n-1}$ is called \emph{admissible} if
\[B_{w_iw_{i+1}}=1, \qquad 0\leq i\leq n-2.\]

\noindent For simplicity, the restriction $\displaystyle \sigma|_{\Sigma_B^+}$ of the left shift to $\Sigma_B^+$ is denoted by $\sigma_B$. The dynamical system $(\Sigma_B^+,\sigma_B)$ is a one-sided \emph{shift of finite type}. In this one-step presentation, admissibility is determined by the finitely many allowed and forbidden transitions encoded by $B$. In general, a shift of finite type is a shift space that can be described by forbidding finitely many finite words.

\medskip
\noindent
The system $(\Sigma_B^+,\sigma_B)$ is called \emph{topologically mixing} if for every pair of nonempty open sets $U,V\subseteq\Sigma_B^+$, there exists $N\geq1$ such that
\[
\sigma_B^n(U)\cap V\neq\varnothing
\qquad
\text{for every }n\geq N.
\]
For a shift of finite type presented by $B$, topological mixing is
equivalent to $B$ being \emph{primitive}; that is, there exists $K\geq1$ such that $B^K>0$.

\noindent Equivalently, for every pair of symbols $i,j\in\mathcal A$, there is an admissible path of length $K$ from $i$ to $j$. This uniform
connecting property allows admissible finite words to be closed into
periodic orbits with a bounded additional transition length. In particular, sufficiently long admissible words can be connected uniformly. This provides a bounded closing length for admissible words
and therefore gives a linear bound on periodic penetration.\\

\noindent We next show that the linear periodic penetration exhibited by the full shift persists for mixing shifts of finite type.

\begin{proposition} [Mixing SFT]
Let $(\Sigma_B^+,\sigma_B)$ be a one-sided mixing shift of finite type,
where $B$ is a primitive transition matrix, and let $\mathcal C_1$
denote the standard one-symbol cylinder partition. Then
$\mathcal C_1$ is $\eta$-admissible. More precisely, there exists
$K\geq1$ such that
\[ M_{\sigma_B}(\mathcal C_n) \leq n+K-1 \quad \text{for every} \quad n\geq 1.\]
\end{proposition}

\begin{proof}
Let $n\in \mathbb{N}$. Since $B$ is primitive, by definition, there exists $K\geq1$ such that
$B^K>0$. Let $w=w_0w_1\cdots w_{n-1}$ be any admissible word of length $n$. Since $B^K>0$, there exists an admissible path of length $K$ from $w_{n-1}$ to $w_0$. Thus, there exist $a_1,\ldots,a_{K-1}\in \mathcal{A}_q$ such that
\[B_{w_{n-1}a_1} = B_{a_1a_2} = \cdots = B_{a_{K-1}w_0} = 1.\]
Consequently, the block $w_0w_1\cdots w_{n-1}a_1\cdots a_{K-1}$ can be repeated periodically to obtain a point $p\in[w]$. Notice that the period of $p$ is at most $n+K-1$. Therefore,
\[\eta_{\sigma_B}([w]) \leq n+K-1.\]

\noindent Because this holds for every nonempty cylinder $[w]\in\mathcal {C}_n$, we obtain $\displaystyle M_{\sigma_B}(\mathcal {C}_n) \leq n+K-1$. Moreover, we have $\displaystyle (\mathcal {C}_1)_{\sigma_B}^{(n)} = \mathcal C_n$.
Hence,
\[ \limsup_{n\to\infty}
\frac{M_{\sigma_B}\bigl((\mathcal {C}_1)_{\sigma_B}^{(n)}\bigr)}{n} =  \limsup_{n\to\infty}
\frac{M_{\sigma_B}\bigl(\mathcal{C}_n\bigr)}{n} \leq
\lim_{n\to\infty} \frac{n+K-1}{n} = 1 <\infty. \]

\noindent Thus, $\mathcal C_1$ is $\eta$-admissible: $\mathcal C_1\in\mathscr A_{\sigma_B}$.
\end{proof}

\begin{remark}
The preceding proposition shows that the linear periodic penetration
observed for the full shift persists for mixing shifts of finite type.
Indeed, every admissible word of length $n$ can be closed into a
periodic orbit by adding at most $K-1$ symbols, where $K$ is independent
of the word and of $n$. Consequently, the maximal least period needed
to penetrate the $n$-cylinders grows at most linearly with $n$.
\end{remark}

\subsection{A Zero-Entropy Example with Dense Periodic Points}\label{subsection: A Zero-Entropy Example with Dense Periodic Points}

We now give an example showing that strong periodic penetration may
occur even when the topological entropy vanishes.

\noindent For each $d\geq1$, let $\displaystyle x^{(d)} = (10^{d - 1})^\infty \in\Sigma_2^+$,
and define \[X = \overline{ \bigcup_{d\geq1}\mathcal O(x^{(d)})}, \qquad
f=\sigma|_X. \]
Here, $\displaystyle (10^{d - 1})^\infty$ means that there are exactly $d-1$ zeros between two consecutive 1's, and $\mathcal O(x^{(d)})$ denotes the forward orbit of $x^{(d)}$ under the shift. Notice that each $x^{(d)}$ has least period $d$.

\begin{proposition}
For the system $(X,f)$ defined above, the set of periodic points is dense in $X$ and $h_{\mathrm{top}}(f)=0$.
Moreover, the standard one-symbol cylinder partition $\mathcal C_1[X] = \mathcal C_1:=\left\{[0], [1]\right\}$ is $\eta$-admissible and satisfies
\[ M_f(\mathcal C_n[X])=n \quad \text{for every} \quad n\geq1,\quad \text{and} \quad A_f(\mathcal C_n[X])\asymp n.\]
Consequently,
\[L_A(f,\mathcal C_1[X])=\infty \quad \text{and hence}, \qquad \Lambda_A(f)=\infty.\]
Here, $\displaystyle \mathcal{C}_n[X]:=\left\{[w]: w\in \{0,1\}^n, [w]\cap X \neq \emptyset \right\}$, and $A_f(\mathcal C_n[X])\asymp n$ means there exist constants $c,C>0$ and $N\in \mathbb{N}$ such that ${\displaystyle cn \leq A_f(\mathcal{C}_n) \leq Cn}$ for every $n\geq N$.
\end{proposition}

\begin{proof}
By definition, $\displaystyle X = \overline{ \bigcup_{d\geq1}\mathcal O(x^{(d)})}$, and every point in $\mathcal O(x^{(d)})$ is periodic. Therefore, $\displaystyle \overline{\operatorname{Per}(f)}=X$. First, we estimate the number of admissible words of length $n$. Let $P(n)=|\mathcal C_n[X]|$. An $n$-word containing at least two occurrences of the symbol $1$ and arising from $\mathcal O(x^{(d)})$, the orbit of $x^{(d)}$, must satisfy $d\leq n-1$. For each such $d$, $\mathcal O(x^{(d)})$ contributes at most $d$ distinct words of length $n$. In addition, there are at most $n$ words containing exactly one occurrence of $1$, together with the all-zero word. Hence, $\displaystyle P(n) \leq 1+n+\sum_{d=1}^{n-1}d$. Thus, $\displaystyle P(n)=O(n^2)$, and consequently $\displaystyle \lim_{n\to\infty}\frac{1}{n}\log P(n)=0$. Since $\mathcal C_1$ is a generating partition, this implies $\displaystyle h_{\mathrm{top}}(f)=0$.

 We next determine the maximal periodic penetration time of the $n$-cylinders. We say the word $w$ is admissible if it occurs in the subshift $X$. Equivalently, the word $w$ is admissible if $\displaystyle [w]\cap X \neq \emptyset$.

If an admissible $n$-word contains at least two occurrences of $1$, then it occurs in the orbit of some $x^{(d)}$ with $d\leq n-1$.
Hence, the corresponding cylinder contains a periodic point of period
at most $n-1$. If an admissible word contains exactly one occurrence of $1$, say in position $i$, where $0\leq i\leq n-1$, then the corresponding cylinder contains a point from the orbit of $x^{(d)}$, where $\displaystyle d=\max\{i+1,n-i\}\leq n$. Finally, the all-zero cylinder contains the fixed point $0^\infty$. Therefore, $\displaystyle M_f(\mathcal C_n[X])\leq n$.

On the other hand, consider the cylinder $[10^{n-1}]$. It contains $x^{(n)}$, which has least period $n$. Moreover, no periodic point of $X$ having least period less than $n$ can belong to this cylinder. Hence
$\eta_f([10^{n-1}])=n$. It follows that $M_f(\mathcal C_n[X])=n$.
Therefore,
\[\limsup_{n\to\infty}
\frac{M_f(\mathcal C_n[X])}{n}=1.\]
And since $\displaystyle (\mathcal C_1)_f^{(n)}=\mathcal C_n[X]$, the partition $\mathcal C_1$ is $\eta$-admissible.

It remains to estimate the average periodic penetration time.
For 
\[\left\lceil\frac{n}{2}\right\rceil \leq d\leq n-1 \qquad \text{and} \qquad  0\leq i\leq n-d-1,\] consider the word
\[w_{d,i} = 0^i1\,0^{d-1}1\,0^{n-i-d-1}.\]
This word contains exactly two occurrences of $1$, separated by
$d$ positions. It occurs in the orbit of $x^{(d)}$, and no periodic
point in $X$ of smaller least period can contain this word. This implies,
$\displaystyle \eta_f([w_{d,i}])=d$.
For each fixed $d$, there are $n-d$ such words. Therefore,
\[\sum_{C\in\mathcal C_n[X]}\eta_f(C) \geq \sum_{d=\lceil n/2\rceil}^{n-1}d(n-d).\]
Because
\[\sum_{d=\lceil n/2\rceil}^{n-1}d(n-d) \asymp n^3 \qquad \text{and} \qquad \left|\mathcal C_n[X]\right|=P(n)=O(n^2),\]
there exists a constant $c>0$ such that, for all sufficiently large
$n$,
\begin{align}\label{A_f(C)>=cn}
    A_f(\mathcal C_n[X])\geq cn.
\end{align}
On the other hand, we have
\begin{align}\label{A_f(C)<=n}
    A_f(\mathcal C_n[X]) \leq M_f(\mathcal C_n[X]) =n
\end{align}
Therefore, by combining (\ref{A_f(C)>=cn}) and (\ref{A_f(C)<=n}), we get
\[A_f(\mathcal C_n[X])\asymp n.\]

\noindent Because $\mathcal C_n[X]$ is a partition, $\displaystyle N(\mathcal C_n[X])=|\mathcal C_n[X]|=P(n)=O(n^2)$. Thus, we have
\[1+\log N(\mathcal C_n[X])=O(\log n) \qquad \text{with} \qquad A_f(\mathcal C_n[X])\geq cn.\]
Hence,
\[L_A(f,\mathcal C_1) = \limsup_{n\to\infty} \frac{A_f(\mathcal C_n[X])}
{1+\log N(\mathcal C_n[X])} = \infty.\]
Since $\mathcal C_1\in\mathscr A_f$, this implies $\Lambda_A(f)=\infty$.
\end{proof}

\begin{remark}
This example demonstrates that periodic penetration and topological
entropy may exhibit substantially different asymptotic behavior.
Although $h_{\mathrm{top}}(f)=0$, the average periodic penetration of the dynamically refined cylinder partitions grows linearly, while their covering complexity grows only polynomially. Consequently, $\Lambda_A(f)=\infty$. Thus, $\Lambda_A$ detects periodic-orbit information that is not captured by topological entropy alone.
\end{remark}

\section{The Invertible Setting}\label{Topological Invertible Setting}

In this section, we consider the case in which $f:X\to X$ is a homeomorphism of a compact metric space $(X,d)$. In this setting, the use of forward joins (Definition \ref{Dynamical Join}) introduces a fundamental restriction that motivated us to adopt and modify our construction.

\noindent Recall that $f$ is called \emph{positively expansive} if there exists $c>0$, for every pair of distinct points $x,y\in X$, there is an integer $n\geq 0$ such that 
\[ d(f^n(x),f^n(y))>c. \]
On the other hand, a homeomorphism $f$ is called \emph{expansive} if the same condition holds with $n\in\mathbb Z$, allowing both forward and backward iterates.

\noindent A classical result due to Schwartzman (see his doctoral dissertation at Yale University in 1952 \cite{schwartzman1952transformation}) states that if $X$ is an infinite compact metric space and $f:X\to X$ is a homeomorphism, then for every $\varepsilon>0$, there exist distinct points
$x,y\in X$ such that
\[\sup_{n\geq 0}d(f^n(x),f^n(y))\leq\varepsilon.\]

\noindent In particular, a homeomorphism of an infinite compact metric space cannot be positively expansive. Equivalently, as stated by Richeson and Wiseman (\cite[Theorem 1]{Positively_expansive_homeomorphisms_of_compact_spaces}), a homeomorphism of a compact metric space is positively expansive if and only if the space is finite.

This observation has a direct consequence for our framework. In fact, the existence of a finite forward generating open cover implies positive expansivity. To see this; suppose that $\mathcal U$ is such a cover, and let $\delta>0$ be a Lebesgue number of
$\mathcal U$. If two points $x,y\in X$ satisfy
\[d(f^n(x),f^n(y))<\delta \qquad\text{for all }n\geq0, \]
then $x$ and $y$ belong to a common element of $\mathcal U_f^{(k)}$ for every $k\geq1$. Since $\mathcal U$ is generating, its forward joins eventually refine open covers of arbitrarily small mesh, forcing $x=y$. Hence $f$ is positively expansive.

\noindent Consequently, whenever $X$ is infinite and $f$ is a homeomorphism, we get
\[\operatorname{Gen}(X,f)=\varnothing, \quad \text{and therefore} \quad \mathscr A_f=\varnothing. \]

\noindent In particular, this applies to the two-sided full shift $(\Sigma_q,\sigma_q)$ for $q\geq2$. Although $\sigma_q$ is expansive, it is not positively expansive and therefore admits no finite forward generating open cover. In fact, the two-sided full shift was the original reason and motivation for this modification.

\noindent The appropriate notion for invertible systems is instead based on two-sided dynamical refinements. In their work on topological entropy and expansiveness, Keynes and Robertson \cite[Theorem~3.2]{Generators_for_topological_entropy_and_expansiveness_Keynes_Robertson_1969} established the connection between expansivity and the existence of finite topological generators for compact dynamical systems. In particular, for a homeomorphism on a compact metric space, expansivity is equivalent to the existence of a finite generating open cover when both positive and negative iterates are allowed.

\noindent Motivated by these results, we replace forward joins with \emph{symmetric joins} defined by
   \[\widehat{\mathcal {U}}_f^{(n)} := \bigvee_{j=-n}^{n}f^{-j}\left(\mathcal U\right), \qquad \text{for} \quad n\geq0. \]
These refinements incorporate information from both the past and the future, allowing us to study infinite expansive homeomorphisms, such as the two-sided full shift.

For finite dynamical systems consisting of periodic orbits, however, the forward and symmetric constructions yield the same invariant. Indeed, every generating cover eventually refines to the singleton partition, and the corresponding penetration times are determined entirely by the least periods of the points. As we shall prove below (Proposition \ref{proposition for Finite periodic systems in the invertible setting}), we have
\[ \widehat{\Lambda}_A(f)=\Lambda_A(f) \] 
for every finite periodic system.

The distinction between forward and symmetric joins is also important when $f$ is not invertible. For a general continuous map $f:X\to X$, the forward joins 
\[ \mathcal U_f^{(n)} = \bigvee_{j=0}^{n-1}f^{-j}\mathcal U \] 
are well-defined, since they involve only inverse images under nonnegative iterates of $f$. However, symmetric joins require both positive and negative iterates of the map and are therefore naturally defined for homeomorphisms. In particular, when $f$ is not invertible, the expressions $f^j$ for $j<0$ do not define single-valued iterates on $X$. Thus, the symmetric construction used here does not apply directly.

\noindent Consequently, we retain the forward-join construction for general continuous maps and introduce symmetric joins specifically for homeomorphisms. The two constructions yield the same invariant for finite periodic systems, as shown in Proposition~\ref{proposition for Finite periodic systems in the invertible setting}. However, no general equality between the resulting invariants is asserted.

\noindent One possible approach to relating the two settings is through the natural extension of a noninvertible dynamical system, which yields a homeomorphism on an inverse-limit space when the original map is surjective. Investigating the relationship between
the corresponding periodic penetration invariants under this construction is beyond the scope of this paper.

We now introduce the symmetric analogues of our generating covers, admissibility condition, and periodic penetration invariant. The preceding construction is based on forward dynamical refinements and is therefore naturally adapted to noninvertible systems. For a homeomorphism, however, both forward and backward orbit information is available. In this setting, it is natural to replace the forward iterated join by a symmetric join that records orbit information in both time directions.

\noindent Throughout this section, let $X$ be a non-empty compact Hausdorff space and let $f:X\to X$ be a homeomorphism.

\subsection{Symmetric Generating Covers and Periodic Penetration}

\begin{definition}[Symmetric Iterated Join]\label{Symmetric Iterated Join}
Let $\mathcal U$ be a finite open cover of $X$. For $n\geq0$, define the \emph{symmetric iterated join} of $\mathcal U$ by
\[\widehat{\mathcal U}_{f}^{(n)} := \bigvee_{j=-n}^{n}f^{-j}\mathcal U. \]
Thus, a typical member of $\widehat{\mathcal U}_{f}^{(n)}$ has the form $\displaystyle \bigcap_{j=-n}^{n}f^{-j}(U_j)$, where $\displaystyle U_j\in\mathcal U$, provided that this intersection is nonempty.
\end{definition}

\begin{lemma}
For every $n\geq0$, we have $\displaystyle \widehat{\mathcal U}_{f}^{(n+1)} \prec \widehat{\mathcal U}_{f}^{(n)}$.
\end{lemma}

\begin{proof}
Since $\displaystyle \{-n,\ldots,n\} \subseteq \{-(n+1),\ldots,n+1\}$, every member of $\widehat{\mathcal U}_{f}^{(n+1)}$ is contained in a member of $\widehat{\mathcal U}_{f}^{(n)}$. Hence, $\displaystyle \widehat{\mathcal U}_{f}^{(n+1)} \prec \widehat{\mathcal U}_{f}^{(n)}$.
\end{proof}

\begin{definition}[Symmetric generating cover]
A finite open cover $\mathcal U$ of $X$ is called a \emph{symmetric generating cover} for $(X,f)$ if, for every finite open cover $\mathcal V$ of $X$, there exists $n\geq0$ such that $\displaystyle \widehat{\mathcal U}_{f}^{(n)}
\prec \mathcal V$.

 We denote the collection of all finite symmetric generating covers by
\[
\widehat{\operatorname{Gen}}(X,f).
\]
\end{definition}

\begin{definition}[Symmetric $\eta$-admissibility]
Let $\displaystyle \mathcal{U}\in \widehat{\operatorname{Gen}}(X,f)$.
We say that $\mathcal U$ is \emph{symmetrically $\eta$-admissible} if
\[ \limsup_{n\to\infty}\frac{M_f\bigl(\widehat{\mathcal U}_{f}^{(n)}\bigr)}{2n+1} <\infty.\]
We denote the collection of all symmetrically $\eta$-admissible
generating covers by
\[
\widehat{\mathscr A}_f
:=
\left\{
\mathcal U\in
\widehat{\operatorname{Gen}}(X,f): \limsup_{n\to\infty}
\frac{M_f\bigl(\widehat{\mathcal U}_{f}^{(n)}\bigr)}{2n+1} <\infty \right\}. \]
\end{definition}

\noindent Notice that this symmetric admissibility condition requires periodic orbit penetration to keep pace, at most linearly, with the length of the orbit window observed by the symmetric iterated joins. Indeed, $\widehat{\mathcal U}_{f}^{(n)}$ records orbit information over the time interval
\[ -n,-n+1,\ldots,0,\ldots,n-1,n, \]
which has length $2n+1$. Thus, symmetric $\eta$-admissibility requires that the largest least period needed to penetrate every member of $\widehat{\mathcal U}_{f}^{(n)}$ to grow at most linearly with the length of this orbit window.

\begin{definition}
Let $\mathcal U\in \widehat{\mathscr A}_f$. We define

\[\widehat{L}_A(f,\mathcal U):=\limsup_{n\to\infty}
\frac{A_f\bigl(\widehat{\mathcal U}_{f}^{(n)}\bigr)}{1+\log N\bigl(\widehat{\mathcal U}_{f}^{(n)}\bigr)}.\]

\noindent The corresponding symmetric periodic penetration invariant is defined by 
\[ \widehat{\Lambda}_A(f) := \sup_{\mathcal U\in \widehat{\mathscr A}_f} \widehat{L}_A(f,\mathcal U). \]
\end{definition}

\noindent
For a fixed symmetrically $\eta$-admissible generating cover
$\mathcal U$, the quantity $\widehat{L}_A(f,\mathcal U)$ measures the asymptotic average periodic penetration of the symmetric dynamical refinements, normalized by their topological complexity. The invariant $\widehat{\Lambda}_A(f)$ takes the supremum of this quantity over all symmetrically $\eta$-admissible generating covers.

\subsection{Conjugacy Invariance}

\begin{theorem}
Let $f:X\to X$ and $g:Y\to Y$ be homeomorphisms, and suppose that $h:X\to Y$ is a topological conjugacy. Then,
\[ \widehat{\Lambda}_A(f) = \widehat{\Lambda_A}(g). \]
\end{theorem}

\begin{proof}
Let $\mathcal U$ be a finite open cover of $X$. By the definition of conjugacy, we have
\[ h\circ f^j=g^j\circ h \qquad \text{for every }j\in\mathbb Z.\]
Thus, $\displaystyle h\bigl(f^{-j}(U)\bigr) = g^{-j}\bigl(h(U)\bigr)$ for every $U\in\mathcal U$ and $j\in\mathbb Z$. Consequently, we have
\[h\bigl(\widehat{\mathcal {U}}_f^{(n)}\bigr) = \widehat{h(\mathcal {U})}_g^{(n)} \qquad \text{for every } n\geq 0.\]

\noindent Since $h$ preserves refinements of open covers,
\[ \mathcal U\in \widehat{\operatorname{Gen}}(X,f)
\quad \text{if and only if} \quad h(\mathcal{U})\in \widehat{\operatorname{Gen}}(Y,g).
\]
Moreover, topological conjugacy preserves least periods of periodic points. Hence,
\[ M_f\bigl(\widehat{\mathcal {U}}_f^{(n)}\bigr) = M_g\bigl(\widehat{h(\mathcal U)}_{g}^{(n)}\bigr),\]
and therefore
\[\mathcal U\in \widehat{\mathscr A}_f \quad \text{if and only if} \quad h(\mathcal U)\in \widehat{\mathscr A}_g. \]

\noindent Similarly, we have 
\[A_f\bigl(\widehat{\mathcal {U}}_f^{(n)}\bigr) = A_g\bigl(\widehat{h(\mathcal U)}_{g}^{(n)}\bigr). \]

\noindent Because $h$ induces a bijection between corresponding members of the two covers and preserves inclusions, we also have
\[N\bigl(\widehat{\mathcal {U}}_f^{(n)}\bigr) = N\bigl(\widehat{h(\mathcal U)}_{g}^{(n)}\bigr).\]
Thus,
\[\widehat{L}_A(f,\mathcal U) = \limsup_{n\to\infty}
\frac{A_f\bigl(\widehat{\mathcal U}_{f}^{(n)}\bigr)}{1+\log N\bigl(\widehat{\mathcal U}_{f}^{(n)}\bigr)}
= \limsup_{n\to\infty}
\frac{A_g\bigl(\widehat{h(\mathcal U)}_{g}^{(n)}\bigr)}{1+\log N\bigl(\widehat{h(\mathcal U)}_{g}^{(n)}\bigr)} =\widehat{L}_A(g,h(\mathcal U)).\]

\noindent Finally, notice that the correspondence $\displaystyle \mathcal U\longmapsto h(\mathcal U)$ is a bijection from $\widehat{\mathscr A}_f$ onto $\widehat{\mathscr A}_g$. By taking suprema, we then obtain
\[\widehat{\Lambda}_A(f) = \widehat{\Lambda}_A(g). \]
\end{proof}

\begin{proposition}[Finite periodic systems in the invertible setting]\label{proposition for Finite periodic systems in the invertible setting}
Let $X$ be a finite union of $r$ distinct periodic orbits with least periods $m_1,\ldots,m_r$, and let $f:X\to X$ be the corresponding dynamical system. Then, $f$ is a homeomorphism, and
\[\widehat{\Lambda}_A(f) = \frac{\displaystyle\sum_{i=1}^{r}m_i^2}
{\displaystyle \left(\sum_{i=1}^{r}m_i\right) \left(1+\log\left(\sum_{i=1}^{r}m_i\right)\right)}.\]

\noindent In particular, $\displaystyle \widehat{\Lambda}_A(f)=\Lambda_A(f)$.
\end{proposition}

\begin{proof}
Since $X$ is finite and Hausdorff, it has the discrete topology. Moreover, $f$ permutes the points of $X$. Thus, it is a homeomorphism.

\noindent Let $\displaystyle \mathcal{U}\in \widehat{\operatorname{Gen}}(X,f)$ be any finite symmetric generating cover of $X$. By the generating property, there exists $N\geq0$ such that $\displaystyle \widehat{\mathcal{U}}_{f}^{(N)}
\prec \bigl\{\{x\}:x\in X\bigr\}$. Since the symmetric joins become finer as $n$ increases, for every $n\geq N$, their distinct nonempty members are precisely the singletons of $X$. For each $x\in X$, we have $\displaystyle \eta_f(\{x\})=\operatorname{per}_f(x)$. This implies
\[M_f\left(\widehat{\mathcal U}_f^{(n)} \right) = \max_{1\leq i\leq r}m_i \qquad \text{for all } n\geq N.\]

\noindent Thus, every finite symmetric generating cover is
$\eta$-admissible. Furthermore, we have
\[N\left(\widehat{\mathcal U}_{f}^{(n)} \right) = |X| = \sum_{i=1}^{r}m_i,\]
and
\[A_f\left(\widehat{\mathcal U}_{f}^{(n)}\right) = \frac{1}{|X|}\sum_{x\in X}\operatorname{per}_f(x) = \frac{\displaystyle \sum_{i=1}^{r}m_i^2}{\displaystyle \sum_{i=1}^{r}m_i}.\]
Therefore,
\[\widehat{L}_A \left(f,\mathcal U \right) = \frac{\displaystyle\sum_{i=1}^{r}m_i^2} {\displaystyle
\left(\sum_{i=1}^{r}m_i\right) \left(1+\log\left(\sum_{i=1}^{r}m_i\right)\right)}.\]

\noindent Since this value is independent of the symmetric generating cover $\mathcal U$, taking the supremum over all admissible symmetric generating covers gives the desired formula. The equality $\displaystyle \widehat{\Lambda}_A(f)=\Lambda_A(f)$ follows from the corresponding formula for finite periodic systems in the forward setting.
\end{proof}

Now we investigate some calculations for the example that motivated us to consider symmetric joins: the 2-sided full shift.

\subsection{The Two-Sided Full Shift}

For $q \geq 2$, let $\displaystyle \Sigma_q = \{0,1,\ldots,q-1\}^{\mathbb Z}$. Let $\displaystyle \sigma_q:\Sigma_q\longrightarrow\Sigma_q$
be the two-sided full $q$-shift (or the two-sided full shift on $q$ symbols). Let $\mathcal C_0$ denote the partition according to the zeroth coordinate:
\[\mathcal C_0 := \left\{\{x\in\Sigma_q:x_0=a\}:a =0, \ldots, q-1\right\}.\]

\noindent The symmetric iterated join
\[\widehat{(\mathcal C_0)}_{\sigma_q}^{(n)} = \bigvee_{j=-n}^{n}\sigma_q^{-j} \left( \mathcal {C}_0 \right) \]
is precisely the partition into cylinders determined by the coordinates $\displaystyle \{-n,-n+1,\ldots,n\}$. Thus, its members can be naturally identified with words of length $2n+1$ over the alphabet $\{0,1,\ldots,q-1\}$.

\begin{theorem}[Two-Sided Full Shifts]\label{Theorem: L_A(sigma,C_0)=1/log q}
For the two-sided full shift $(\Sigma_q,\sigma_q)$, the partition
$\mathcal C_0$ is a symmetric generating partition and is
symmetrically $\eta$-admissible. Moreover,
\[M_{\sigma_q}\left(\widehat{(\mathcal C_0)}_{\sigma_q}^{(n)}\right) = 2n+1 \qquad \text{and} \qquad \displaystyle \lim_{n \rightarrow \infty }\frac{A_{\sigma_q}
\left(\widehat{(\mathcal C_0)}_{\sigma_q}^{(n)}\right)}{2n+1} = 1.\]
Consequently, $\displaystyle \widehat{L}_A(\sigma_q,\mathcal C_0) = \frac{1}{\log q}$, and therefore
\[ \widehat{\Lambda}_A(\sigma_q) \geq \frac{1}{\log q}. \]
\end{theorem}

\begin{proof}
Since the symmetric iterated joins of $\mathcal C_0$ determine
arbitrarily large finite coordinate windows, $\mathcal C_0$ is a
symmetric generating partition. Let $n\geq 0$ be fixed, and set $m=2n+1$. Notice that each member of $\displaystyle \widehat{(\mathcal C_0)}_{\sigma_q}^{(n)}$ corresponds to a word
$\displaystyle w=w_0w_1\cdots w_{m-1}$ of length $m$. Let $p(w)$ denote the least positive integer $d\leq m$ such that
\[ w_i=w_{i+d} \qquad \text{whenever }0\leq i<m-d.\]
Thus, $p(w)$ is the least period of the finite word $w$. We claim that
\[\eta_{\sigma_q}([w])=p(w),\]
where $[w]$ denotes the corresponding cylinder. Indeed, if $w$ has
period $d$, then repeating its first $d$ symbols produces a bi-infinite periodic point of period dividing $d$ contained in $[w]$. Conversely, if $[w]$ contains a periodic point of least period $d$, then the coordinates of $w$ satisfy $w_i=w_{i+d}$
whenever both sides occur in the prescribed coordinate window. Hence, $p(w)\leq d$. This implies $\displaystyle \eta_{\sigma_q}([w])=p(w)$.

Notice that every word of length $m$ has period at most $m$, and the word $\displaystyle 0^{m-1}1$ has least period exactly $m$. Therefore,
\[M_{\sigma_q} \left(\widehat{(\mathcal C_0)}_{\sigma_q}^{(n)}
\right) = m = 2n+1.\]
Consequently, we have
\[\limsup_{n\to\infty}\frac{M_{\sigma_q}\left(\widehat{(\mathcal C_0)}_{\sigma_q}^{(n)}\right)}{2n+1} = 1.\]
So, $\mathcal C_0$ is symmetrically $\eta$-admissible:  $\displaystyle \mathcal{C}_0 \in \widehat{\mathscr A}_{\sigma_q}$.

Now we estimate the average periodic penetration time. Since the symmetric join contains exactly $q^m$ cylinders, we have
\begin{align}\label{A_sigma_q}
    A_{\sigma_q}\left(\widehat{(\mathcal C_0)}_{\sigma_q}^{(n)}
\right) = \frac{1}{q^m} \sum_{w\in\{0,\ldots,q-1\}^m}p(w).
\end{align}
For each $1\leq d<m$, there are at most $q^d$ words of length $m$
having least period $d$, since such a word is completely determined by its first $d$ symbols. 

\noindent Therefore,
\begin{align*}
0 &\leq m- A_{\sigma_q}\left(\widehat{(\mathcal C_0)}_{\sigma_q}^{(n)}\right)\\
&= m- \frac{1}{q^m} \sum_{w\in\{0,\ldots,q-1\}^m}p(w), \quad \text{from } (\ref{A_sigma_q})\\
&=\frac{1}{q^m}\sum_{w\in\{0,\ldots,q-1\}^m}\bigl(m-p(w)\bigr)\\
&\leq \frac{1}{q^m} \sum_{d=1}^{m-1}(m-d)q^d.
\end{align*}

\noindent Let $k=m-d$. Then,
\[\frac{1}{q^m}\sum_{d=1}^{m-1}(m-d)q^d = \sum_{k=1}^{m-1}\frac{k}{q^k} \leq \sum_{k=1}^{\infty}\frac{k}{q^k} = \frac{q}{(q-1)^2}. \]
Hence,
\[m-\frac{q}{(q-1)^2} \leq A_{\sigma_q} \left( \widehat{(\mathcal C_0)}_{\sigma_q}^{(n)} \right) \leq m. \]
By dividing $m=2n+1$ and taking the limit, we then obtain
\[ \lim_{n\rightarrow \infty} \frac{A_{\sigma_q}\left(\widehat{
(\mathcal C_0)}_{\sigma_q}^{(n)}\right)}{2n+1} = 1.\]

\noindent Because the symmetric join is a partition into exactly
$q^{2n+1}$ cylinders, we have 
\[N\left(\widehat{(\mathcal C_0)}_{\sigma_q}^{(n)}\right) = q^{2n+1}.\]

\noindent Therefore,
\begin{align}
\nonumber \widehat{L}_A(\sigma_q,\mathcal C_0) &= \limsup_{n\to\infty}
\frac{A_{\sigma_q} \left(\widehat{(\mathcal C_0)}_{\sigma_q}^{(n)}\right)}{1+\log N\left( \widehat{(\mathcal C_0)}_{\sigma_q}^{(n)}\right)}\\
\nonumber &= \lim_{n\to\infty} \frac{(2n+1)+O(1)} {1+(2n+1)\log q}\\
&= \frac{1}{\log q}.\label{L_hat(sigma_q, C_1)}
\end{align}
Because $\displaystyle \mathcal C_0\in \widehat{\mathscr A}_{\sigma_q}$ it follows from the definition of the supremum that
\[ \widehat{\Lambda}_A(\sigma_q) \geq \frac{1}{\log q}.\]
\end{proof}

\begin{remark}
The preceding calculation shows that the symmetric formulation achieves, for the two-sided full shift, the same canonical value obtained from the standard cylinder generator of the one-sided full shift. From (\ref{L_A(sigma_q, C_1)}) and (\ref{L_hat(sigma_q, C_1)}) we see that
\[  L_A(\sigma_q,\mathcal C_1) = \frac{1}{\log q} \qquad \text{and} \qquad \widehat{L}_A(\sigma_q,\mathcal C_0) = \frac{1}{\log q}. \]
The symmetric join is essential in the 2-sided shift because it records both past and future coordinates of the invertible shift.
\end{remark}

\section{Variance of Local Periodic Penetration Times}\label{sec:variance}

The invariant $\Lambda_A(f)$ measures the average scale of the least periods needed to encounter the elements of increasingly refined covers. The average alone, however, does not record how unevenly these least periods are distributed among the cover elements. We, therefore, introduce a variance-based companion invariant. Throughout this section, we retain the standing assumptions and notation of the preceding sections, in particular the class $\mathscr A_f$ of admissible forward generating covers and the quantities $\eta_f(U)$, $A_f(\mathcal U)$, and $N(\mathcal U)$.

\begin{definition}[Periodic Penetration Variance]\label{Definition Penetration-Variance}
For a finite open cover $\mathcal U$ of $X$, define its \emph{periodic penetration variance} by
\[ V_f(\mathcal U) :=\frac{1}{|\mathcal U|}\sum_{U\in\mathcal U}
 \bigl(\eta_f(U)-A_f(\mathcal U)\bigr)^2.\]
 
\noindent For $\mathcal U\in\mathscr A_f$, let
\[ L_V(f,\mathcal U) :=\limsup_{n\to\infty} \frac{V_f(\mathcal U_f^{(n)})} {\bigl(1+\log N(\mathcal U_f^{(n)})\bigr)^2}
 \qquad \text{and} \qquad \Lambda_V(f):=\sup_{\mathcal U\in\mathscr A_f}L_V(f,\mathcal U).\]
\end{definition}

\noindent The squared denominator puts the variance on the square of the logarithmic covering-complexity scale used in $\Lambda_A$. In particular, $\Lambda_V$ records the dispersion of local least periods, not the number of periodic points in individual open sets. As with $\Lambda_A$, the supremum is taken over all admissible generators; computing $L_V$ for a particular generator does not, by itself, determine $\Lambda_V$.

\begin{proposition}[Conjugacy invariance]\label{proposition variance conjugacy}
Let $X$ and $Y$ be nonempty compact Hausdorff spaces, and let $f:X\to X$ and $g:Y\to Y$ be continuous maps with dense periodic points. Suppose that $(X,f)$ and $(Y,g)$ are topologically
conjugate and that $\mathscr A_f\neq\varnothing$. Then,
\[\Lambda_V(f)=\Lambda_V(g).\]
\end{proposition}

\begin{proof}
Let $h:X\to Y$ be a conjugacy, so $h\circ f=g\circ h$.
For each nonempty open $U\subseteq X$, conjugacy preserves the least period of every periodic point in $U$. Hence, $\eta_g(h(U))=\eta_f(U)$. Thus, for any finite open cover $\mathcal U$, the covers $\mathcal U$ and $h(\mathcal U):=\{h(U):U\in\mathcal U\}$ have the same penetration times, counted with multiplicity. Consequently,
\[ A_g(h(\mathcal U))=A_f(\mathcal U) \qquad \text{and} \qquad
 V_g(h(\mathcal U))=V_f(\mathcal U).\]

\noindent Moreover, we have $h(\mathcal U_f^{(n)})=(h(\mathcal U))_g^{(n)}$ for all $n\geq 0$, and $N(h(\mathcal W))=N(\mathcal W)$ and $M_g(h(\mathcal W))=M_f(\mathcal W)$ for every finite open
cover $\mathcal W$. The generating and admissibility conditions are therefore preserved by $h$. Hence, $L_V(f,\mathcal U)=L_V(g,h(\mathcal U))$ for every $\mathcal U\in\mathscr A_f$. By taking suprema and applying the same argument to $h^{-1}$, we obtain the claim.
\end{proof}

\subsection{Finite periodic systems}

\begin{proposition}\label{proposition variance-finite}
Suppose $X$ is a finite discrete space and $f:X\to X$ is a permutation whose distinct cycles have lengths $m_1,\ldots,m_r$. Set $\displaystyle M=\sum_{i=1}^r m_i$. Then,
\[\displaystyle  \Lambda_V(f)=  \frac{1}{(1+\log M)^2}
 \left( \frac{\sum_{i=1}^r m_i^3}{M} -\left(\frac{\sum_{i=1}^r m_i^2}{M}\right)^2 \right).\]
In particular, $\Lambda_V(f)=0$ if and only if all the cycle lengths $m_1,\ldots,m_r$ are equal.
\end{proposition}
\begin{proof}
Every forward generating cover of a finite discrete system has an
iterated join equal to the singleton partition. One can see this by applying the generating property to the singleton cover and using finiteness. Subsequent joins remain the singleton partition because $f$ is a permutation. All such generators are admissible. Notice that, at a singleton belonging to a cycle of length $m_i$, the penetration time is $m_i$; this value occurs $m_i$ times. Therefore, the mean and second moment of the penetration
times of the singleton partition are $\displaystyle  \frac{\sum_i m_i^2}{M}$ and $\displaystyle \frac{\sum_i m_i^3}{M}$ respectively. Subtracting the square of the mean and dividing by $(1+\log M)^2$ gives us the formula. The variance vanishes precisely when all the singleton penetration times, equivalently all the cycle lengths, are equal.
\end{proof}

\begin{remark}
For two cycles of lengths $m$ and $n$, the formula simplifies to
\[ \Lambda_V(f)=  \frac{mn(m-n)^2}{(m+n)^2(1+\log(m+n))^2}.\]

\noindent Thus, the variance invariant detects unequal cycle lengths even though every finite periodic system has zero topological entropy.
\end{remark} 

\subsection{Symbolic examples}

\begin{proposition}[Canonical cylinders of the full shift]
\label{Proposition: variance fullshift}
Let $\displaystyle \sigma_q:\Sigma_q^+\to\Sigma_q^+$ be the one-sided full shift on $q\ge2$ symbols, and let $\mathcal C_1$ be its one-symbol cylinder partition. Then,
\[ \sup_{n\ge1}V_{\sigma_q}(\mathcal C_n)<\infty \qquad \text{and} \qquad L_V(\sigma_q,\mathcal C_1)=0.\]
\end{proposition}

\begin{proof}
Let $\mathcal C_n$ denote the partition of the one-sided full shift $\Sigma_q^+$ into cylinders of length $n$. For a word $w=w_1\cdots w_n$, a \emph{proper border} of length $\ell$, where $1\leq \ell<n$, is a prefix of $w$ that is also a suffix of $w$. Equivalently, $\displaystyle w_i=w_{n-\ell+i}$, for $1\leq i\leq\ell$. This condition is equivalent to saying that $w$
has combinatorial period $n-\ell$.

\noindent Define
\[D_n(w):=n-\eta_{\sigma_q}([w]).\]
Since the least period of a periodic point in $[w]$ is the least combinatorial period of $w$, the quantity $D_n(w)$ is precisely the length of the longest proper border of $w$. We set $D_n(w)=0$
when $w$ has no proper border.

\noindent Equip the finite set of words $\{0,\ldots,q-1\}^n$ with the uniform probability measure. Thus, each word has probability $q^{-n}$. We write $\mathbb P$ for probability and $\mathbb E$ for expectation with respect to this measure. For $1\leq\ell<n$, let $B_\ell$ be the event that $w$ has a border of length $\ell$. As mentioned above, the border condition is equivalent to
\[w_i=w_{i+n-\ell}, \qquad 1\leq i\leq\ell.\]
Once the first $n-\ell$ symbols are specified, all remaining symbols are uniquely determined. In other words, exactly $q^{n-\ell}$ words satisfy this condition. Hence,
\[ \mathbb P(B_\ell) = \frac{q^{n-\ell}}{q^n}=q^{-\ell}.\]

\noindent For $1\leq k<n$, the event $\{D_n\geq k\}$ is the union of the events $B_\ell$ over $k\leq\ell<n$. By the union bound, we get 
\begin{align*}
\mathbb P(D_n\geq k)
&\leq\sum_{\ell=k}^{n-1}\mathbb P(B_\ell)\\
&=\sum_{\ell=k}^{n-1}q^{-\ell}\\
&\leq\sum_{\ell=k}^{\infty}q^{-\ell}
=\frac{q^{1-k}}{q-1}.
\end{align*}

\noindent Next, we estimate the second moment of $D_n$; that is, the expected value of $\displaystyle D_n^2$. Notice that for every nonnegative integer $k$, we have $\displaystyle k^2=\sum_{m=1}^{k}(2m-1)$.

\noindent Let $\displaystyle \Omega_n:=\{0,\ldots,q-1\}^n$, and for each $k\in\{0,\ldots,n-1\}$, let $\displaystyle A_k:=\{w\in\Omega_n:D_n(w)=k\}$.

\noindent Since every word in $\Omega_n$ has probability $q^{-n}$,
we obtain
\begin{align*}
\mathbb E[D_n^2]
&= \sum_{w\in\Omega_n}D_n(w)^2\mathbb P(\{w\})\\
&= \frac{1}{q^n} \sum_{k=0}^{n-1}\sum_{w\in A_k}D_n(w)^2\\
&= \frac{1}{q^n}\sum_{k=0}^{n-1}\sum_{w\in A_k}k^2\\
&= \sum_{k=0}^{n-1}k^2\frac{|A_k|}{q^n}\\
&= \sum_{k=0}^{n-1}k^2\mathbb P(D_n=k)\\
&= \sum_{k=1}^{n-1}\left(\sum_{m=1}^{k}(2m-1)\right)
   \mathbb P(D_n=k)\\
&= \sum_{m=1}^{n-1}(2m-1)\sum_{k=m}^{n-1}\mathbb P(D_n=k)\\
&= \sum_{m=1}^{n-1}(2m-1) \mathbb P(D_n\geq m)\\
&\leq \frac{q}{q-1}\sum_{m=1}^{n-1}(2m-1)q^{-m}\\
&\leq \frac{q}{q-1}\sum_{m=1}^{\infty}(2m-1)q^{-m}\\
&= \frac{q(q+1)}{(q-1)^3}.
\end{align*}

\noindent In particular, the second moments of $D_n$, that is $\displaystyle \mathbb{E}[D_n^2]$, are bounded uniformly in $n$. Since $\displaystyle \eta_{\sigma_q}([w])=n-D_n(w)$, and the cylinders in $\mathcal C_n$ correspond bijectively to the $q^n$ words of length $n$, we have
\[V_{\sigma_q}(\mathcal C_n) =\operatorname{Var}(D_n) = \mathbb{E}[D_n^2] - \left(\mathbb{E}[D_n]\right)^2\leq\mathbb{E}[D_n^2] \leq\frac{q(q+1)}{(q-1)^3}. \]

\noindent Finally, the one-symbol cylinder cover
$\mathcal C_1$ is generating and admissible, and
\[\mathcal C_{1,\sigma_q}^{(n)} = \mathcal C_n, \qquad N(\mathcal C_n)=q^n.\]
It follows that
\[0 \leq \lim_{n\rightarrow \infty}\frac{V_{\sigma_q}(\mathcal C_n)}{(1+\log N(\mathcal C_n))^2} \leq \lim_{n\rightarrow \infty}\frac{q(q+1)}{(q-1)^3(1+n\log q)^2}= 0.\]
Hence,
\[L_V(\sigma_q,\mathcal C_1)=0.\]
\end{proof}

\begin{remark}
Proposition \ref{Proposition: variance fullshift} computes the value for the canonical cylinder generator only. It does \emph{not} establish that $\Lambda_V(\sigma_q)=0$, since the definition of $\Lambda_V$ involves a supremum over all admissible open generating covers.
\end{remark}

\begin{proposition}[A zero-entropy subshift with divergent variance] \label{proposition variance sparse}
For $d\ge1$, let $x^{(d)}=(10^{d-1})^\infty$. Let
\[ X=\overline{\bigcup_{d\ge1}\mathcal O(x^{(d)})}
 \subseteq\{0,1\}^{\mathbb N_0},\qquad f=\sigma|_X.\]
 
\noindent Then, $X$ has dense periodic points and zero topological entropy. For its partition $\mathcal C_n(X)$ into nonempty length-$n$ cylinders, we have $\displaystyle |\mathcal C_n(X)|=\left\lfloor\frac{n^2}{4}\right\rfloor+n+1$. This implies
\[ A_f(\mathcal C_n(X)) =\frac n2+O(1) \qquad \text{and}\qquad
 V_f(\mathcal C_n(X)) =\frac{n^2}{24}+O(n).\]
Consequently,
\[ L_V(f,\mathcal C_1(X))=\Lambda_V(f)=\infty.\]
\end{proposition}

\begin{proof}
The periodic points are dense by the definition of $X$. Every allowed word of length $n$ is either all zeros, has exactly one
symbol $1$, or has at least two symbols $1$ occurring at a fixed 
spacing $d\in\{1,\ldots,n-1\}$. In the last case, the word is
determined by its spacing $d$ and the residue class of the positions of its $1$'s modulo $d$. Exactly $\min\{d,n-d\}$ residue classes produce at least two occurrences. Thus, with
$P_n:=|\mathcal C_n(X)|$, we get
\[ P_n=1+n+\sum_{d=1}^{n-1}\min\{d,n-d\} = 1+n+\left\lfloor\frac{n^2}{4}\right\rfloor = \frac{n^2}{4}+O(n).\]

\noindent A word with at least two $1$'s at spacing $d$ belongs to a periodic orbit of least period $d$, and no periodic point of smaller period can contain that word. Its cylinder therefore has penetration time $d$. A word containing exactly one $1$, at position $i\in\{0,\ldots,n-1\}$, has penetration time $\max\{i+1,n-i\}$, which is at most $n$, and the all-zero word has penetration time $1$.

\noindent By using the Riemann-sums, we have the estimation
\[ \sum_{d=1}^{n-1}d^k\min\{d,n-d\} =n^{k+2}\int_0^1 t^k\min\{t,1-t\}\,dt +O(n^{k+1}) \qquad \text{for }k=0,1,2.\]
The three integrals, for $k=0,1,2$, are
$1/4$, $1/8$, and $7/96$, respectively.

\noindent The $n$-number of single-$1$ cylinders contribute at most $O(n^2)$ to the sum of penetration times and at most $O(n^3)$ to the sum of their squares. Therefore, we have
 \begin{align*}
 \sum_{C\in\mathcal C_n(X)}\eta_f(C)
 &=\frac{n^3}{8}+O(n^2), \qquad{and}\\
 \sum_{C\in\mathcal C_n(X)}\eta_f(C)^2
 &=\frac{7n^4}{96}+O(n^3).
 \end{align*}
 
\noindent By dividing by $P_n=n^2/4+O(n)$, we get
\[ A_f(\mathcal C_n(X))=\frac n2+O(1) \qquad \text{and} \qquad
 \frac1{P_n}\sum_{C\in\mathcal C_n(X)}\eta_f(C)^2
 =\frac{7n^2}{24}+O(n).\]
Hence,
\[ V_f(\mathcal C_n(X)) =\left(\frac7{24}-\frac14\right)n^2+O(n)
 =\frac{n^2}{24}+O(n).\]
 
\noindent Because the collection of cylinder sets form a partition, $N(\mathcal C_n(X))=P_n$, and therefore
\[\lim_{n\rightarrow \infty} \frac{V_f(\mathcal C_n(X))}
 {(1+\log N(\mathcal C_n(X)))^2}  = \lim_{n\rightarrow \infty}\frac{n^2}{96(\log n)^2} = \infty.\]
 
\noindent Notice that the one-symbol cylinder partition is a forward generator. It is admissible because $M_f(\mathcal C_n(X))=n$ for every $n\ge1$; all penetration times are at most $n$, and the cylinder $[10^{n-1}]\cap X$ has penetration time $n$. Thus, we have $\Lambda_V(f)=\infty$. Finally, $P_n=O(n^2)$ implies 
\[h_{\mathrm{top}}(f)=\lim_{n\to\infty}\frac{1}{n}\log P_n=0.\]
\end{proof}

\begin{remark}
The sparse subshift also satisfies $L_A(f,\mathcal{C}_1(X))=\infty$, since its mean penetration time
is asymptotic to $\displaystyle \frac{n}{2}$, whereas $\displaystyle \log P_n\sim2\log n$. The two examples above show that the normalized variance for a canonical cylinder generator can vanish or diverge. They do not resolve the possible finite positive values of $\Lambda_V$ for infinite subshifts.
\end{remark}

\subsection{The Invertible Setting}

When $f:X\to X$ is a homeomorphism, the symmetric construction from the preceding section has a direct variance analogue. For a symmetric admissible generator $\mathcal U\in \widehat{\mathscr A}_f$, define
\[ \widehat{L}_V(f,\mathcal U)  :=\limsup_{n\to\infty}
 \frac{V_f( \widehat{\mathcal U_{f}}^{(n)})}
 {\bigl[1+\log N(\widehat{\mathcal U_f}^{(n)})\bigr]^2}
 \qquad \text{and} \qquad  \widehat{\Lambda}_V(f) := \sup_{\mathcal U\in \widehat{\mathscr A}_f}
 \widehat{L}_V(f,\mathcal U).\]
 
\noindent We use similar arguments as in the proof of Proposition \ref{proposition variance conjugacy} for symmetric joins and show that $\widehat{\Lambda}_V$ is a conjugacy invariant.

\begin{proposition}\label{proposition variance symmetric finite}
If $X$ is finite and $f:X\to X$ is a permutation, then
\[ \widehat{\Lambda}_V(f)=\Lambda_V(f).\]
For the two-sided full $q$-shift with its zero-coordinate cylinder
partition $\mathcal C_0$, one has
\[\widehat{L}_V(\sigma_q,\mathcal C_0)=0.\]
\end{proposition}
\begin{proof}
For a finite permutation, every symmetric generating cover eventually refines the singleton partition, so the proof of
Proposition \ref{proposition variance-finite} applies unchanged.
For the two-sided full shift, the $n$-th symmetric join of $\mathcal C_0$ is the partition into coordinate blocks indexed by
$\{-n,\ldots,n\}$. There are $q^{2n+1}$ number of such blocks, and the minimum period among all periodic sequences containing
a given block equals the least period of the corresponding finite word of length $2n+1$. Thus, the variance agrees with that of the length-$(2n+1)$ cylinder partition in the one-sided full shift. Proposition \ref{Proposition: variance fullshift} bounds this
variance uniformly, whereas the squared normalizing denominator is
$(1+(2n+1)\log q)^2$, which diverges to $\infty$ as $n \rightarrow \infty$.
\end{proof}

\begin{remark}
The examples illustrate the different roles of the two invariants $\Lambda_A$ and $\Lambda_V$. The invariant $\Lambda_A$ measures normalized average periodic penetration times, while the invariant $\Lambda_V$ measures their normalized heterogeneity. In particular, vanishing variance for one generating cover must not be confused with the vanishing of the supremum defining $\Lambda_V$.
\end{remark}

\begin{remark}[Weighted periodic penetration averages]
The average periodic penetration time can be generalized
by assigning different weights to the elements of a
finite open cover. For $\alpha\in\mathbb R$, define
\[
A_f^{(\alpha)}(\mathcal U)
:=
\frac{\displaystyle\sum_{U\in\mathcal U}
\eta_f(U)^{\alpha+1}}
{\displaystyle\sum_{U\in\mathcal U}
\eta_f(U)^\alpha}.
\]
In particular, $A_f^{(0)}(\mathcal U)=A_f(\mathcal U)$.

These weighted averages are preserved under topological
conjugacy. Indeed, if $h:X\to Y$ is a topological
conjugacy between $f$ and $g$, then
$\eta_f(U)=\eta_g(h(U))$ for every $U\in\mathcal U$,
and hence
\[
A_f^{(\alpha)}(\mathcal U)
=
A_g^{(\alpha)}(h(\mathcal U)).
\]

\noindent Moreover, for every $\alpha\in\mathbb R$,
\[\min_{U\in\mathcal U}\eta_f(U)
\le A_f^{(\alpha)}(\mathcal U)
\le M_f(\mathcal U),\]
and the function
$\alpha\mapsto A_f^{(\alpha)}(\mathcal U)$
is non-decreasing.

\noindent Of particular interest is the case $\alpha=1$. Using the definition of $V_f(\mathcal U)$, we obtain

\[A_f^{(1)}(\mathcal U) = \frac{\displaystyle\sum_{U\in\mathcal{U}} \eta_f(U)^2}
{\displaystyle\sum_{U\in\mathcal{U}}\eta_f(U)} = \frac{V_f(\mathcal{U})+A_f(\mathcal{U})^2}{A_f(\mathcal{U})} = A_f(\mathcal{U}) +\frac{V_f(\mathcal{U})}{A_f(\mathcal{U})}.\]

\noindent Thus, the weighted average corresponding to $\alpha=1$
is determined by the ordinary average and variance.
More generally, weighted averages emphasize different
parts of the distribution of periodic penetration times.
\end{remark}

\section{Measure-Theoretic Return-Time Averages}
\label{section: Measure-Theoretic Return-Time Averages}

We now consider a measure-theoretic analogue of the preceding
constructions. Since individual periodic points may have measure zero, and since periodic-point data are not generally preserved under measure-theoretic isomorphisms, we replace periodic penetration times by positive-measure overlap times.

Throughout this section, $(X,\mathcal B,\mu,f)$ denotes a
probability preserving system. Thus, $f:X\to X$ is measurable and
$\mu(f^{-1}B)=\mu(B)$ for every $B\in\mathcal B$.
All measurable partitions are considered modulo null sets, and
null atoms are discarded. All logarithms are natural.

\subsection{Motivation: Recurrence, Observation, and Information}
\label{subsection:measure-motivation}

The purpose of this construction is to measure recurrence of
\emph{observable sets}, with respect to the measure, rather than recurrence of individual points. In a topological system, a nonempty open set may contain periodic points whose least periods provide a measure of how deeply periodic orbits penetrate the space. In a measure-preserving system, however, individual periodic points may have measure zero. Thus, even the presence of periodic points may not, by itself, invariant under measure-theoretic isomorphism. The natural replacement is to ask whether a set overlaps its own inverse image with \emph{positive
measure}.\\

\noindent There are three different notions that must not be confused:
\begin{enumerate}
    \item [(i).] the least lag (time step or iteration) at which some positive-measure collection of points returns to a set,
    \item[(ii).]  the first-return time of an individual point, and
    \item[(iii).]  the expected first-return time of a
random point in that set.
\end{enumerate}

Our construction uses the first of these. In particular, it measures an \emph{earliest feasible recurrence} and not a typical waiting time.

For a finite measurable partition $\displaystyle \mathcal{P}=\{P_1, \ldots, P_r\}$ of $X$, its atoms represent finitely many observable outcomes. The iterated join $\displaystyle \mathcal P_f^{(n)}=\bigvee_{j=0}^{n-1}f^{-j}\mathcal P$ distinguishes points according to their first $n$ observations under the dynamics. More precisely, each atom of $\displaystyle \mathcal{P}^{(n)}_f$ consists of points sharing the same finite itinerary through the partition elements. Such itineraries may be represented by words of length $n$ over the alphabet $\{1, \dots, r\}$, regardless of whether the underlying dynamical system is symbolic. The positive-measure atoms may have different minimum positive-measure return times. This variation motivates studying not only their measure-weighted average but also their variance and maximum.

\begin{definition}[Positive measure overlap return time]
\label{def:measure-return}
For $B\in\mathcal B$ with $\mu(B)>0$, let
\[
\tau_\mu(B):=\min\{k\ge1:\mu(B\cap f^{-k}B)>0\}.
\]
\end{definition}

Poincar\'e recurrence ensures (See Corollary \ref{overlap finiteneww}) that $\tau_\mu(B)<\infty$.
This quantity is the \emph{least lag (iteration) admitting a positive-measure self-overlap}. It is not the expected first-return time appearing in Kac's formula (see Theorem \ref{Kac's theorem} below).

\subsection{Classical recurrence and the distinction from first returns}
\label{subsection:classical-recurrence}

We recall the precise recurrence statements needed to interpret
$\tau_\mu$. For $B\in\mathcal B$ with $\mu(B)>0$, let
\[
 r_B(x):=\inf\{n\ge1:f^n(x)\in B\},\qquad x\in B,
\]
with the convention that the infimum is $+\infty$ if no return
occurs. Thus, $r_B$ is a \emph{pointwise first-return time}, whereas $\tau_\mu(B)$ is a single number attached to the entire set.

\begin{theorem}[Poincar\'e recurrence theorem]
\label{Poincare Recurrence Theorem}
If $f$ preserves a finite measure $\mu$, then for every measurable
$B$ with $\mu(B)>0$, almost every $x\in B$ returns to $B$
infinitely many times. In particular, $r_B(x)<\infty$ for
$\mu$-almost every $x\in B$.
\end{theorem}

\begin{corollary}\label{overlap finiteneww}
For every $B\in\mathcal B$ of positive measure,
$1\le\tau_\mu(B)<\infty$.
\end{corollary}
\begin{proof}
By Poincar\'e recurrence, the set of points in $B$ that return
at least once has measure $\mu(B)>0$. Since
\[
 \{x\in B:r_B(x)<\infty\}
 =\bigcup_{k\ge1}(B\cap f^{-k}B),
\]
countable additivity implies that at least one intersection on
the right has positive measure. The minimum of the corresponding
positive integers is therefore finite.
\end{proof}

\begin{theorem}[Kac's return-time formula] \label{Kac's theorem}
Suppose $f$ is an ergodic probability-preserving transformation
and $\mu(B)>0$. Then,
\[
 \int_B r_B(x)\,d\mu(x)=1
 \qquad \text{and} \qquad
 \mathbb E_{\mu(\cdot\mid B)}[r_B]
 =\frac1{\mu(B)}.
\]
Here, $\displaystyle \mathbb E_{\mu(\cdot\mid B)}[r_B]$ is the expected value of the first return time $r_B$, and $\displaystyle \mu (\cdot \mid B)$ is the normalized restriction of $\mu$ to $B$. So, $\displaystyle \mathbb E_{\mu(\cdot\mid B)}[r_B] := \int_B r_B(x)\ d\mu_B(x) = \frac{1}{\mu(B)} \int_B r_B(x)\ d\mu(x)$.
\end{theorem}

\noindent The ergodicity hypothesis is important for this formulation.
For a general probability-preserving system, the return-time
integral equals the measure of the invariant saturation of $B$; $\displaystyle Sat_f(B):=\bigcup_{n=0}^\infty f^{-n}(B)$, (the smallest $f$-invariant subset containing $B$, up to null sets) which need not be $1$.

\begin{proposition}[Minimum Overlap Vs Mean First Return]
\label{prop:overlap-versus-kac}
For every positive-measure set $B$, we have
\[
 \tau_\mu(B)=\operatorname*{ess\,inf}_{x\in B}r_B(x),
\]
where the essential infimum is with respect to the conditional
measure $\mu(\cdot\mid B)$. Consequently,
\[
 \tau_\mu(B)\le
 \mathbb E_{\mu(\cdot\mid B)}[r_B]
\]
whenever the latter is defined. In the ergodic case,
$\tau_\mu(B)\le1/\mu(B)$.
\end{proposition}

\begin{proof}

Let $m=\tau_\mu(B)$. If $m=1$, then clearly we have $\displaystyle \tau_\mu(B)=\text{ess} \inf_{x\in B} r_B(x)$ . So, suppose $m\geq 2$. Then, $\displaystyle \mu \left( B\cap f^{-m}(B) \right)>0$ and $\displaystyle \mu \left( B\cap f^{-j}(B) \right)=0$ for all $j=1, \ldots, m-1$. This implies $\displaystyle \mu \left( \{x\in B: r_B(x)=m\} \right)>0$ and $\displaystyle \mu \left( \{x\in B: r_B(x)=j\} \right)=0$ for all $j=1, \ldots, m-1$. Because $\displaystyle B\cap f^{-m}(B) \subseteq \bigcup_{j=1}^m \{x\in B: r_B(x)=j\}$, we then have 
\begin{align*}
    B\cap f^{-m}(B) \overset{\mu-\text{a.e}}{\subseteq} \{x\in B: r_B(x)=m\}.
\end{align*}
On the other hand,
\begin{align*}
    \{x\in B: r_B(x)=m\} \subseteq B\cap f^{-m}(B).
\end{align*}
Thus, $\displaystyle \mu \left( \left(B\cap f^{-m}(B) \triangle \{x\in B: r_B(x)=m\} \right) \right)=0$, and therefore $\displaystyle  \tau_\mu(B)=\operatorname*{ess\,inf}_{x\in B}r_B(x)$.
The expectation bound follows immediately from the definition.
\end{proof}

\begin{example}[A sharp difference from Kac's mean]
\label{ex:kac-versus-overlap}
In the fair binary Bernoulli shift, let $B=[0^n]$. Then
$\mu(B)=2^{-n}$, while $B\cap\sigma^{-1}B$ has positive
measure, so $\tau_\mu(B)=1$. Kac's formula gives
\[
 \mathbb E_{\mu(\cdot\mid B)}[r_B]=2^n.
\]
Thus, the earliest possible positive-measure return occurs after
one step, even though the mean first-return time is exponential in
$n$. This gives an example of why we cannot describe the present invariant as a normalized mean waiting time.
\end{example}

\noindent Let \(\mathcal P=\{P_1,\ldots,P_r\}\) be a finite measurable partition of \(X\). Each atom represents a measurable collection of points whose recurrence can be studied through its minimum positive-measure return time. Different atoms may have different return times, and their measures indicate the relative sizes of the corresponding collections of points. It is therefore natural to investigate both the typical magnitude of these return times and the extent to which they vary among the partition atoms. We introduce the following quantities to capture these features.

\begin{definition}[Mean, Variance, and Maximum]
\label{def:measure-average}
For a finite measurable partition $\mathcal P$, define
\begin{align*}
 A_\mu(\mathcal P)
  &:=\sum_{P\in\mathcal P}\mu(P)\tau_\mu(P),\\
 V_\mu(\mathcal P)
  &:=\sum_{P\in\mathcal P}\mu(P)
       \bigl(\tau_\mu(P)-A_\mu(\mathcal P)\bigr)^2,\\
 M_\mu(\mathcal P)
  &:=\max_{P\in\mathcal P}\tau_\mu(P), \quad \text{and}\\
 H_\mu(\mathcal P)
  &:=-\sum_{P\in\mathcal P}\mu(P)\ln\mu(P).
\end{align*}
\end{definition}

\subsection{Interpretation of the partition statistics}
\label{subsection:measure-interpretation}

The number $A_\mu(\mathcal P)$ is the probability-weighted
average of the earliest feasible return lags (time steps, or iterations) of the partition elements (one may use the terms atoms or observation classes). It gives greater importance to frequently observed
atoms. The number $M_\mu(\mathcal P)$ is a worst-case statistic:
it is sensitive to even one extremely rare positive-measure atom
with a long minimal overlap time. Finally,
$V_\mu(\mathcal P)$ quantifies the dispersion of these lags.
A vanishing variance means that all positive-measure atoms have
the same minimum overlap time; a small \emph{normalized}
variance has the weaker meaning that fluctuations are small
relative to the chosen information scale.

\begin{example}[Zero and positive variance of minimum return times]
Consider the identity map $f=\operatorname{id}_X$ on a probability space $(X,\mathcal B,\mu)$. Every measurable set $B$ with $\mu(B)>0$ satisfies $\displaystyle \tau_\mu(B)=1$. Consequently, for every finite measurable partition
$\mathcal P$, we have
\[A_\mu(\mathcal P)=M_\mu(\mathcal P)=1, \qquad \text{and} \qquad
V_\mu(\mathcal P)=0.\]

\noindent By contrast, consider a finite probability space $X=\{1,2,3,4,5\}$ equipped with the uniform
probability measure, and let $p:X\to X$ be the
permutation defined by
\[p: =(1\ 2)(3\ 4\ 5).\]
For the partition into singletons $\displaystyle \mathcal P=\bigl\{\{1\},\{2\},\{3\},\{4\},\{5\}\bigr\}$, the minimum return times are
\[\tau_\mu(\{1\})=\tau_\mu(\{2\})=2 \qquad \text{and} \qquad
\tau_\mu(\{3\})=\tau_\mu(\{4\}) =\tau_\mu(\{5\})=3. \]
Thus, we have
\[A_\mu(\mathcal P) = \frac25\cdot 2+\frac35\cdot 3 =\frac{13}{5},
\]
and
\begin{align*}
V_\mu(\mathcal P) =\frac25\left(2-\frac{13}{5}\right)^2
 +\frac35\left(3-\frac{13}{5}\right)^2 = \frac{6}{25}>0.
\end{align*}
Moreover, $M_\mu(\mathcal P)=3$.

\noindent These examples illustrate that the variance $V_\mu(\mathcal P)$ measures the variability of minimum positive-measure return times among the partition elements, whereas $A_\mu(\mathcal P)$ records only their probability weighted average.
\end{example}

\begin{proposition}\label{prop:measure-bounds}
If $\mathcal P$ has $r$ number of positive measure atoms, then
\[A_\mu(\mathcal P)\le r \qquad \text{and} \qquad 1\le A_\mu(\mathcal P)\le M_\mu(\mathcal P)\le
\frac1{\min_{P\in\mathcal P}\mu(P)}.\]
\end{proposition}
\begin{proof}
Let $P\in \mathcal{P}$ be fixed, and let $m=\tau_\mu(P)$. For $0\le i<j<m$, measure preservation gives
\[\mu(f^{-i}(P)\cap f^{-j}(P)) =\mu(P\cap f^{-(j-i)}(P))=0.\] 
 This implies that the sets $P,f^{-1}(P),\ldots,f^{-(m-1)}(P)$ are pairwise disjoint modulo null sets. Hence, $\mu(P)\tau_\mu(P) = m\mu(P)\le1$. By summing
$\mu(P)\tau_\mu(P)\le1$ over the $r$ atoms, we
get $A_\mu(\mathcal P)\le r$. The other inequalities follow directly from the definitions and the fact $1\le\tau_\mu(P)\le1/\mu(P)$.
\end{proof}

\begin{proposition}[Isomorphism covariance]
\label{propposition measure isomorphism}
Suppose $(X,\mathcal B,\mu,f)$ and $(Y,\mathcal C,\nu,g)$ are
measure-theoretically isomorphic. Let $\phi:X \rightarrow Y$ be a measure theoretic isomorphism (modulo null sets), and $\mathcal P$ is a finite measurable partition,
then
\[\tau_\mu(B)=\tau_\nu(\phi B),\qquad
(A_\mu,V_\mu,M_\mu,H_\mu)(\mathcal P)
=(A_\nu,V_\nu,M_\nu,H_\nu)(\phi\mathcal P).\]
\end{proposition}
\begin{proof}
Let $k\ge1$ and let $B\in \mathcal{B}$. Then,
\begin{align*}
    \mu(B\cap f^{-k}(B)) &= \mu \left[B\cap \phi^{-1}\circ g^{-k}(\phi (B)) \right]\\
    &= \mu \left[(\phi^{-1} \circ \phi (B)) \cap (\phi^{-1}\circ g^{-k}(\phi (B))) \right], \because\ \phi\ \text{is a bijection}\\
    &= \mu \left[\phi^{-1} (\phi (B) \cap  g^{-k}(\phi (B))) \right]\\
    &= \nu \left[\phi (B) \cap  g^{-k}(\phi (B)) \right],\ \because\ \phi\ \text{is a measure isomorphism}.
\end{align*}

So, 
\[
\mu(B\cap f^{-k}(B)) =\nu(\phi (B) \cap g^{-k}\circ \phi (B)) \qquad \text{for every } B\in \mathcal{B}. \]
Thus, corresponding atoms have identical overlap times and masses.
Therefore, the four asserted identities follow immediately from their definitions.
\end{proof}

\subsection{Asymptotic invariants and return-time admissibility}

Let $(X,\mathcal B,\mu,f)$ be a measure preserving probability space, and let $\mathcal P=\{P_1,\ldots,P_r\}$ be a finite measurable partition of $X$. For each $n\geq 1$, define the $n$-fold dynamical refinement of $\mathcal P$ by
\[\mathcal P_f^{(n)} := \bigvee_{j=0}^{n-1}f^{-j}(\mathcal P).\]
Here, $f^{-j}(\mathcal P)$ denotes the partition $\displaystyle f^{-j}(\mathcal P) = \{f^{-j}(P_i):1\leq i\leq r\}$, and the join $\displaystyle \bigvee$ denotes the common refinement
of the indicated partitions. Thus, the atoms of $\mathcal P_f^{(n)}$ are the nonempty sets of the form
\[P_{i_0}\cap f^{-1}(P_{i_1})
\cap\cdots\cap f^{-(n-1)}(P_{i_{n-1}}),\]
where $i_0,\ldots,i_{n-1}\in\{1,\ldots,r\}$.
Such an atom consists of points whose first $n$
iterates visit the partition elements
$P_{i_0},\ldots,P_{i_{n-1}}$, respectively.

\noindent We say that $\mathcal P$ is a \emph{forward generating partition} if
\begin{align}\label{generating partition}
\sigma\left(
\bigcup_{n=1}^{\infty}\mathcal P_f^{(n)}
\right) = \mathcal B
\quad(\mathrm{mod}\ \mu).
\end{align}
Here, $\sigma(\cdot)$ denotes the sigma-algebra
generated by a collection of measurable sets (not the shift map). Equivalently, $\mathcal P$ is forward generating if the sigma-algebra generated by
\[\{f^{-j}(P_i):j\geq 0,\ 1\leq i\leq r\}\]
coincides with $\mathcal B$ modulo $\mu$. Loosely speaking, this means that observing which element of $\mathcal P$ contains $f^j(x)$, for every $j\geq 0$, provides
enough information to distinguish measurable sets up to sets of measure zero.

\noindent As we did in Section \ref{Topological Invertible Setting}, for an invertible probability-preserving
transformation, one may instead use the customary \emph{two-sided generating partition}
condition
\[\sigma\left(
\bigcup_{j\in\mathbb Z}f^{-j}\mathcal P
\right) = \mathcal B \quad(\mathrm{mod}\ \mu). \]
In this case, observations are allowed in both forward and backward time.

\noindent Every forward generating partition is two-sided generating when $f$ is invertible, but the converse need not hold. Thus, the forward-generator condition is stronger for invertible systems. Throughout this paper, unless otherwise specified, the term \textit{generating partition} will refer to a forward generating partition. This convention applies to both invertible and non-invertible probability-preserving systems.

\subsection{Why normalize by Shannon entropy?}
\label{subsection:measure-normalization}

The Shannon entropy $H_\mu(\mathcal P_f^{(n)})$ measures the
average information required to distinguish the positive-measure
observation words of length $n$. For more details, please refer to his original paper \cite{shannon1948mathematical}, and other standard texts such as \cite{books_for_Shannon_entropy_cover2006elements}, \cite{mane1987ergodic}. For a finite partition, the sequence is subadditive, and hence the entropy rate
\[
 h_\mu(f,\mathcal P):=
 \lim_{n\to\infty}\frac1n H_\mu(\mathcal P_f^{(n)})
\]
exists by the F\`ekete's lemma \cite{fekete1923verteilung}. The Kolmogorov-Sinai entropy is $\displaystyle h_\mu(f) :=\sup_{\mathcal P}h_\mu(f,\mathcal P)$, with the supremum being taken over all finite measurable partitions. For more details, one may refer to standard texts such as \cite{kolmogorov1958new}, \cite{sinai1959concept}, \cite{walters_book_1982introduction}. A finite forward generator satisfies $h_\mu(f,\mathcal P)=h_\mu(f)$
by the Kolmogorov-Sinai Generator Theorem. The ratio
\[ \frac{A_\mu(\mathcal P_f^{(n)})}
 {1+H_\mu(\mathcal P_f^{(n)})}\]
therefore compares an \emph{earliest feasible recurrence scale}
with an \emph{average observational-information scale}. Notice that the additive $1$ avoids a zero denominator for trivial partitions and finite zero-entropy systems; it does not, however, change the limit when $H_\mu(\mathcal P_f^{(n)})\to\infty$.

\noindent The variance is divided by the square of the same scale because
variance is quadratic in return times. This choice gives us
\[ \frac{V_\mu(\mathcal P_f^{(n)})}
 {(1+H_\mu(\mathcal P_f^{(n)}))^2}\]
as a dimensionless measure of \emph{relative} dispersion.(How much the values vary or spread around their mean.) It is clear that neither ratio should be interpreted as entropy itself, nor as
a universal recurrence law independent of the chosen partition.

\noindent When the entropy rate is positive and $A_\mu(\mathcal P_f^{(n)})=n+O(1)$, the ratio tends to $1/h_\mu(f,\mathcal P)$. A smaller value may, therefore , indicate \emph{more information per unit observation time}, not necessarily a smaller absolute return time. The Bernoulli calculations below (Theorem \ref{theorem: Canonical Bernoulli calculations}) make this distinction explicit.

\begin{definition}[Linear Maximal Return-Time Admissibility]
\label{def: Return-Time Admissibility}
A finite measurable partition $\mathcal P$ is said to be
\emph{return-time admissible} if it is generating and
$\displaystyle \limsup_{n\to\infty}
\frac{M_\mu(\mathcal {P}_f^{(n)})}{n}<\infty$. Equivalently, $\displaystyle M_\mu(\mathcal {P}_f^{(n)}) = O(n)$.

Let $\mathscr A_\mu(f)$ denote the class of such partitions.
\end{definition}

\begin{remark}[Motivation for Return-Time Admissibility.]
\label{remark: why admissibility?}
The generating condition ensures that the partition captures all the measurable information about the system, up to sets of measure zero (modulo null sets). Thus, we consider only partitions that describe the entire system rather than just part of it.

\noindent The linear maximal-return-time condition requires that there exist a constant $C_{\mathcal P}>0$ such that, for all sufficiently large $n$, every positive measure atom $B$ of $\mathcal P_f^{(n)}$ satisfies $\displaystyle \tau_\mu(B)\le nC_{\mathcal P}$. The constant may depend on the partition $\mathcal P$. This condition is stronger than requiring only
$\displaystyle A_\mu(\mathcal P_f^{(n)})=O(n)$, because the average may remain small even when some atoms with very small probabilities have large return times. 

\noindent We impose this additional condition by choice; it does not follow from the generating property. However, this restriction has some limitations. It may exclude certain finite generating partitions, and some systems may have no finite generating partition satisfying the condition. Furthermore, even if each admissible partition gives a finite value for an invariant, the supremum over all admissible partitions may still be infinite. However, these issues do not affect the proof of measure theoretic isomorphism invariance, provided the admissible class is preserved under isomorphism.
\end{remark}

\begin{definition}[Asymptotic mean and variance invariants]
\label{def:measure-asymptotic}
For a finite partition $\mathcal P$, let
\[L_\mu(f,\mathcal P)
 :=\limsup_{n\to\infty}
 \frac{A_\mu(\mathcal P_f^{(n)})}
 {1+H_\mu(\mathcal P_f^{(n)})} \qquad \text{and} \qquad
L_{V,\mu}(f,\mathcal P)
 :=\limsup_{n\to\infty}
 \frac{V_\mu(\mathcal P_f^{(n)})}
 {\bigl(1+H_\mu(\mathcal P_f^{(n)})\bigr)^2}.\]
 
\noindent If $\mathscr A_\mu(f)\ne\varnothing$, define
\[
\Lambda_\mu(f):=
\sup_{\mathcal P\in\mathscr A_\mu(f)}L_\mu(f,\mathcal P) 
\qquad \text{and} \qquad \Lambda_{V,\mu}(f):=
\sup_{\mathcal P\in\mathscr A_\mu(f)}L_{V,\mu}(f,\mathcal P).\]
\end{definition}
\begin{remark}\label{rem:admissibility-limits}
These global quantities may take the value $\infty$. If $\mathscr A_\mu(f)=\varnothing$, we leave them undefined. 

\noindent The return-time condition implies that, for each admissible $\mathcal{P}$, the maximum overlap time grows at most linearly. It does \emph{not} impose a uniform constant over all admissible partitions. Consequently, it does not, by itself, establish that the global suprema are finite. It is not presently established here whether the condition is satisfied by every finite measurable Bernoulli generator.
\end{remark}

\begin{proposition}[Finite partition-level values at positive entropy] \label{proposition: measure positive entropy}
Suppose that $\mathcal P$ is an admissible finite
generator and that $0<h_\mu(f)<\infty$. Then,
\[L_\mu(f,\mathcal P)<\infty
\qquad\text{and}\qquad L_{V,\mu}(f,\mathcal P)<\infty.\]
\end{proposition}

\begin{proof}
Let $\displaystyle \mathcal Q_n=\mathcal P_f^{(n)}$ and $\displaystyle H_n=H_\mu(\mathcal Q_n)$.
Since $\mathcal P$ is admissible, there exist
constants $C>0$ and $N_1\in \mathbb{N}$ such that $\displaystyle M_\mu(\mathcal Q_n)\le nC$ for all $\displaystyle n\ge N_1$. Recall that, for any finite measurable partition $\mathcal Q$,
the definitions of the average and variance give: 
\[A_\mu(\mathcal Q)\le M_\mu(\mathcal Q) \qquad \text{and} \qquad
V_\mu(\mathcal Q)
\le \sum_{B\in\mathcal Q}
\mu(B)\tau_\mu(B)^2
\le M_\mu(\mathcal Q)^2.
\]
Consequently, for all $n\ge N_1$, $\displaystyle A_\mu(\mathcal Q_n)\le nC$ and $\displaystyle V_\mu(\mathcal Q_n)\le n^2C^2$.

\noindent On the other hand, since $\mathcal P$ is generating, the Kolmogorov-Sinai generator theorem implies $\displaystyle 
\lim_{n\to\infty}\frac{H_n}{n} =h_\mu(f)>0$.
Therefore, there exists $N_2\ge 1$ such that
\[H_n\geq \frac{nh_\mu(f)}{2},
\qquad\text{for all }n\ge N_2.\]

\noindent Hence, for every $n\ge\max\{N_1,N_2\}$, we have 
\[\frac{A_\mu(\mathcal Q_n)}{H_n}
\le \frac{2C}{h_\mu(f)} \qquad \text{and} \qquad
\frac{V_\mu(\mathcal Q_n)}{H_n^2}
\le \frac{4C^2}{h_\mu(f)^2}.\]
Taking the limit superior in both inequalities, we conclude that
\[
L_\mu(f,\mathcal P)<\infty
\qquad\text{and}\qquad
L_{V,\mu}(f,\mathcal P)<\infty.
\]
\end{proof}

\begin{theorem}[Measure-Theoretic Isomorphism Invariance]
\label{Theorem: Measure-Theoretic Isomorphism Invariance}
Let $(X,\mathcal B,\mu,f)$ and $(Y,\mathcal C,\nu,g)$ be measure-theoretically isomorphic probability preserving dynamical systems.
Then, $(X,\mathcal B,\mu,f)$ admits an admissible finite generator if and only if $(Y,\mathcal C,\nu,g)$ does. Furthermore, whenever the global invariants are defined, we have
\[\Lambda_\mu(f)=\Lambda_\nu(g) \qquad\text{and}\qquad \Lambda_{V,\mu}(f)=\Lambda_{V,\nu}(g).\]
\end{theorem}

\begin{proof}
Let $(X,\mathcal B,\mu,f)$ and $(Y,\mathcal C,\nu,g)$ be measure-theoretically isomorphic probability-preserving dynamical systems,
and let $\phi:X_0\to Y_0$ be an isomorphism between invariant subsets of full measure. Let $\mathcal P$ be a finite measurable partition of $X$, and let $\phi (\mathcal P)$ denote the corresponding partition of $Y$, defined modulo null sets. Since $\phi$ is a measure-preserving bijection with a measurable inverse, it preserves the measures of the corresponding partition atoms. Moreover, the conjugacy relation $\displaystyle \phi\circ f=g\circ\phi$ implies,
\[\phi\left(\mathcal P_f^{(n)}\right) = (\phi (\mathcal P))_g^{(n)} \quad(\mathrm{mod}\ \nu) \qquad \text{for every} \qquad n\geq 1. \]

\noindent Because $\phi$ and $\phi^{-1}$ are measurable, the sigma-algebra generated by the dynamical refinements of $\mathcal P$ coincides with $\mathcal B$ (modulo $\mu$) if and only if the
corresponding sigma-algebra generated by the refinements of $\phi (\mathcal P)$ coincides with $\mathcal C$ (modulo $\nu$).
Thus, $\mathcal P$ is a finite generator if and
only if $\phi(\mathcal P)$ is a finite generator.

\noindent By Proposition \ref{propposition measure isomorphism},
for every $n\ge 1$, we have
\begin{align*}
M_\mu(\mathcal P_f^{(n)})
&=M_\nu((\phi (\mathcal P))_g^{(n)}),\\
A_\mu(\mathcal P_f^{(n)})
&=A_\nu((\phi (\mathcal P))_g^{(n)}),\\
V_\mu(\mathcal P_f^{(n)})
&=V_\nu((\phi (\mathcal P))_g^{(n)}), \quad \text{and}\\
H_\mu(\mathcal P_f^{(n)})
&=H_\nu((\phi (\mathcal P))_g^{(n)}).
\end{align*}

\noindent In particular, we have 
\[M_\mu(\mathcal P_f^{(n)})=O(n) \qquad \text{if and only if} \qquad M_\nu((\phi\mathcal P)_g^{(n)})=O(n).\]
Thus, $\mathcal P$ is admissible if and only if $\phi (\mathcal P)$ is admissible. Notice that the equalities above also imply
\[ L_\mu(f,\mathcal P) = L_\nu(g,\phi(\mathcal P)) \qquad \text{and} \qquad L_{V,\mu}(f,\mathcal P) = L_{V,\nu}(g,\phi (\mathcal P)),\]
whenever these quantities are defined. Finally, because $\phi$ establishes a bijective correspondence between the admissible finite generators of the two systems, taking suprema over these classes yields
\[\Lambda_\mu(f)=\Lambda_\nu(g) \qquad \text{and} \qquad
\Lambda_{V,\mu}(f)=\Lambda_{V,\nu}(g),\]
whenever defined. Therefore, the existence of an admissible
finite generator and the values of the corresponding global invariants are preserved under measure-theoretic isomorphism.
\end{proof}

\subsection{Bernoulli factors and the danger of unrestricted suprema}
\label{subsection: Bernoulli factors and the danger of unrestricted suprema}

The following classical theorem helps explain why it is natural
to distinguish generating partitions from arbitrary finite
partitions. We state it in its standard invertible, ergodic
form; no assertion about non-invertible factors is needed here.

\begin{theorem}[Sinai's Factor Theorem \cite{Sinai1964}, for the earliest version of the theorem, see \cite{Sinai1962}]
\label{thm:sinai-factor}
Let $(X,\mathcal B,\mu,f)$ be an ergodic invertible
probability-preserving transformation of a standard probability
space with $h_\mu(f)>0$. For every $0<t<h_\mu(f)$ there is
a nontrivial Bernoulli shift of entropy $t$ that is a
measure-theoretic factor of $(X,\mu,f)$.
\end{theorem}

\noindent Recall that we say $\left(Y, \mathcal{C}, \nu, g \right)$ is a measure theoretic factor of $\left(X, \mathcal{B}, \mu, f \right)$, if there exists a semi-conjugacy $\displaystyle \pi: X \rightarrow Y$.

\begin{remark}
In Sinai's Factor Theorem, an invertible probability-preserving transformation is a measure-preserving map $f:X\to X$ that admits
a measurable inverse modulo sets of measure zero. The invertibility assumption allows the Bernoulli factor to be represented by a two-sided shift. A corresponding version of the theorem also holds for noninvertible ergodic probability preserving transformations, with a one-sided Bernoulli shift
as the factor. In either case, the existence of a Bernoulli factor
means that there is a measurable, measure-preserving
semiconjugacy from the original system onto the
Bernoulli shift. It does not imply that the two
systems are measure theoretically isomorphic (measure theoretically conjugate).
\end{remark}

Let $h>0$ be the entropy of a two-sided Bernoulli shift, and let
$0<t<h$.  By Sinai's Factor Theorem, there exists a factor map onto a Bernoulli shift of entropy $t$. Taking the inverse images of the atoms of its canonical one-symbol partition gives a finite measurable partition $\mathcal {R}_t$, say, of the original system. This partition $\mathcal R_t$ cannot be forward generating.

\noindent Let $\mathcal Q$ denote the canonical one-symbol
partition of the Bernoulli factor, and let $\pi:X\to Y$ be the corresponding factor map. Since $\mathcal R_t=\pi^{-1}(\mathcal Q)$, we have
\[H_\mu((\mathcal R_t)_f^{(n)}) = H_\nu(\mathcal Q_g^{(n)}), \qquad \text{for every } n\geq 1.\] 
Consequently,
\[h_\mu(f,\mathcal R_t) = h_\nu(g,\mathcal Q) = t < h_\mu(f).\]
Thus, by the Kolmogorov-Sinai generator theorem, $\mathcal R_t$ cannot be forward generating.

Since the factor map preserves the measures of corresponding partition atoms and their positive-measure return times, the recurrence statistics of $\mathcal {R}_t$ can be computed using the corresponding Bernoulli partition. The canonical Bernoulli calculation below (See Theorem \ref{theorem: Canonical Bernoulli calculations}) then gives $\displaystyle  L_\mu(f,\mathcal R_t)=\frac{1}{t}$. As $t\searrow 0$, these values diverge. Consequently, if one defines an \emph{unrestricted} global mean invariant by
\[ \Lambda_\mu^{\mathrm{all}}(f):=  \sup_{\mathcal R\text{ finite}}L_\mu(f,\mathcal R),\]
then for every positive-entropy, ergodic invertible system
covered by Sinai's theorem, will have 
\[\Lambda_\mu^{\mathrm{all}}(f)=\infty.\]

\noindent The partitions $\mathcal R_t$ are generally not forward
generating, as established above.
This motivates restricting our attention to forward
generating partitions. However, this argument does
not establish that the supremum over linearly
return-time admissible partitions is finite, nor
does it justify the necessity of the additional
linear maximal-return-time condition: $\displaystyle M_\mu(\mathcal P_f^{(n)})=O(n)$ (Definition \ref{def: Return-Time Admissibility}).

\begin{remark}[Scope of the Factor Argument]
The strict inequality $t<h_\mu(f)$ is sufficient for the preceding argument. We restrict our attention to invertible probability-preserving systems so that the stated version of Sinai's Factor Theorem works.

Although the partitions $\mathcal R_t$ may satisfy
the linear maximal-return-time condition, they
are not forward generating (as shown above). Consequently, they do not belong to the class $\mathscr A_\mu(f)$ of linearly return-time admissible partitions. Thus, the factor argument alone does not determine the supremum over $\mathscr A_\mu(f)$.
\end{remark}

\subsection{Bernoulli shifts}
\label{subsection:measure-bernoulli}

Let $q\ge2$, let $\Sigma_q^+=\{0,\ldots,q-1\}^{\mathbb N_0}$,
and let $\sigma$ be the left shift. Fix a probability vector
$\mathbf p=(p_0,\ldots,p_{q-1})$ with all $p_i>0$, and equip
$\Sigma_q^+$ with the Bernoulli product measure $\mu_{\mathbf p}$.
Let
\[H(\mathbf p):=-\sum_{i=0}^{q-1}p_i\log p_i>0 \qquad \text{and} \qquad \rho:=\max_i p_i<1.\]
Let $\mathcal C=\{[0],\ldots,[q-1]\}$ be the canonical
one-symbol partition, and let $\mathcal C_n=\mathcal C_\sigma^{(n)}$.

\begin{lemma}\label{lemma: Bernoulli Return}
For each word $w$ of length $n$,
\[\tau_{\mu_{\mathbf p}}([w])=\operatorname{per}(w),\]
where $\operatorname{per}(w)$ is the smallest period of $w$. In-other words, 
\[\operatorname{per}(w):= \min\{1\le k\le n: w_j=w_{j+k}\ \text{for all }0\le j<n-k\}.\]

\noindent In particular, we have
\[
M_{\mu_{\mathbf p}}(\mathcal C_n)=n.
\]
\end{lemma}
\begin{proof}
For $k<n$, the two length-$n$ cylinders intersect with
positive measure precisely when their overlapping symbols agree.
For $k=n$, the blocks are disjoint and independence gives
$\mu_{\mathbf p}([w]\cap\sigma^{-n}[w])=
\mu_{\mathbf p}([w])^2>0$.
Thus, the least overlap iteration time is the least combinatorial period. The word $0^{n-1}1$ has least combinatorial period $n$.
\end{proof}

\begin{theorem}[Canonical Bernoulli calculations]
\label{theorem: Canonical Bernoulli calculations}
The canonical partition $\mathcal C$ is admissible. For every
$n\ge1$, we have 
\[n-\frac{\rho}{(1-\rho)^2}
 \le A_{\mu_{\mathbf p}}(\mathcal C_n)\le n \qquad \text{and} \qquad  V_{\mu_{\mathbf p}}(\mathcal C_n)
 \le \frac{\rho(1+\rho)}{(1-\rho)^3}.\]
Consequently,
\[L_{\mu_{\mathbf p}}(\sigma,\mathcal C)
=\frac1{H(\mathbf p)} \qquad \text{and} \qquad L_{V,\mu_{\mathbf p}}(\sigma,\mathcal C)=0,\]
and hence
\[\Lambda_{\mu_{\mathbf p}}(\sigma)\ge\frac1{H(\mathbf p)} \qquad \text{and} \qquad \Lambda_{V,\mu_{\mathbf p}}(\sigma)\ge0.\]
\end{theorem}
\begin{proof}
Let $W_n$ be a random word of length $n$ whose letters are
independent and identically distributed according to the probability distribution $\mathbf p$. Notice that
\begin{align*}
 A_{\mu_{\mathbf p}}(\mathcal C_n) &= \sum_{B\in \mathcal{C}_n} \mu_{\mathbf p}(B)  \tau_{\mu_{\mathbf p}} (B)\\
 &=\sum_{w\in \mathcal{A}^n} \mu_{\mathbf p}([w])  \tau_{\mu_{\mathbf p}} ([w])
\end{align*}

\noindent By Lemma \ref{lemma: Bernoulli Return}, we have $\displaystyle \tau_{\mu_{\mathbf p}}([w]) = \operatorname{per} (w)$. Thus,
\begin{align*}
 A_{\mu_{\mathbf p}}(\mathcal C_n)  =\sum_{w\in \mathcal{A}^n} \mu_{\mathbf p}([w])  \operatorname{per} (w).
\end{align*}
Notice that $W_n$ is a random word of length $n$ whose letters are independently distributed according to the probability vector $\mathbf p$. Thus, $\displaystyle \mathbb{P}(W_n=w) = \mu_{\mathbf p}([w])$. Therefore,
\begin{align*}
    A_{\mu_{\mathbf p}}(\mathcal C_n)  &= \sum_{w\in \mathcal{A}^n} \mu_{\mathbf p}([w])  \operatorname{per} (w)\\
    &=\sum_{w\in \mathcal{A}^n} \mathbb{P}(W_n=w) \operatorname{per} (w)\\
    &= \mathbb{E}[\operatorname{per}(W_n)], \quad \text{by the definition of expectation for a discrete random variable.}
\end{align*}
So,
\begin{align} \label{Average = Expectation}
     A_{\mu_{\mathbf p}}(\mathcal C_n)  = \mathbb E[\operatorname{per}(W_n)].
\end{align}

\noindent We closely follow the arguments used in the proof of Proposition \ref{Proposition: variance fullshift}. Let $D_n:=n-\operatorname{per}(W_n)$. This is the length of the longest proper border of $W_n$, where a border is a
proper prefix that is also a suffix. Therefore,
\begin{align}\label{A(C_n) = n-E[D(n)]}
    A_{\mu_{\mathbf p}}(\mathcal C_n)=n-\mathbb E[D_n].
\end{align}

\noindent For $1\le\ell<n$, let $B_\ell$ be the event that $W_n$
has a border of length $\ell$. Writing $d=n-\ell$, the
border condition is then equivalent to
\[ W_j=W_{j+d}, \quad \text{for } 0\le j<n-d.\]

\noindent Let $\displaystyle \rho := \max_{1\leq i\leq m}p_i<1$. Once the first $d$ symbols are specified, the remaining
$\ell$ symbols are uniquely determined by these relations. Since the letters are independent, the conditional probability that
each subsequent symbol equals its prescribed value, given the preceding symbols, is at most $\rho$. Therefore, $\displaystyle \Pr(B_\ell)\le\rho^\ell$. For $1\le k<n$, the union bound implies
\begin{align}\label{bound for Pr(D_n>=k)}
     \Pr(D_n\geq k)  \leq \sum_{\ell=k}^{n-1}\Pr(B_\ell)  \leq \sum_{\ell=k}^{\infty}\rho^\ell  = \frac{\rho^k}{1-\rho}.
\end{align}

\noindent Since $\displaystyle D_n\in\mathbb{Z}^+_0$ with $D_n\le n-1$, we have
\[D_n=\sum_{k=1}^{n-1}\mathbf 1_{\{D_n\ge k\}}.\]
Here, 
\[
\mathbf{1}_{\{D_n\ge k\}}=
\begin{cases}
1, & D_n\ge k,\\
0, & D_n<k.
\end{cases}
\]

\noindent By taking expectations and using the linearity, we obtain
\[\mathbb E[D_n] =\sum_{k=1}^{n-1}\Pr(D_n\ge k).\]
Therefore, by the preceding estimate, we get
\[\mathbb E[D_n]
\le \frac{1}{1-\rho}\sum_{k=1}^{\infty}\rho^k
=\frac{\rho}{(1-\rho)^2}.\]

\noindent Since $D_n\ge0$, we obtain the two-sided estimate
\[ n-\frac{\rho}{(1-\rho)^2}
 \le A_{\mu_{\mathbf p}}(\mathcal C_n)\le n.\]

\noindent For the variance, we have
\begin{align*}
    V_{\mu_{\mathbf p}}(\mathcal C_n) &:= \sum_{w\in \mathcal{A}^n} \mu_{\textbf{p}} ([w]) \left(\tau_{\mu_{\textbf{p}}} ([w]) - A_{\mu_{\textbf{p}}}(\mathcal{C}_n) \right)^2\\
    &= \sum_{w\in \mathcal{A}^n} \mu_{\textbf{p}} ([w]) \left(\operatorname{per}(w) - A_{\mu_{\textbf{p}}}(\mathcal{C}_n) \right)^2, \quad \text{by Lemma }\ref{lemma: Bernoulli Return}\\
    &= \sum_{w\in \mathcal{A}^n} \mu_{\textbf{p}} ([w]) \left(\operatorname{per}(w) - \mathbb{E}[\operatorname{per} (W_n)]\right)^2, \quad \text{by }(\ref{Average = Expectation})\\
    &=: \operatorname{Var}(\operatorname{per} (W_n)).
\end{align*}

\noindent Recall that $\displaystyle D_n=n-\operatorname{per}(W_n)$. Since subtracting a random variable from a constant does not change its variance, we then get
\begin{align}\label{bound estimate for variance}
    V_{\mu_{\mathbf p}}(\mathcal C_n) = \operatorname{Var}(\operatorname{per} (W_n)) =\operatorname{Var}(D_n) =\mathbb E[D_n^2]-(\mathbb E[D_n])^2 \le \mathbb E[D_n^2].
\end{align}

\noindent To estimate $\mathbb E[D_n^2]$, observe that for
every nonnegative integer $d$, $\displaystyle d^2=\sum_{k=1}^{d}(2k-1)$. Since $D_n\in \mathbb{Z}^+_0$ and $0\le D_n\le n-1$, it follows that
\[D_n^2 =\sum_{k=1}^{n-1}(2k-1) \mathbf 1_{\{D_n\ge k\}}. \]

\noindent Thus,
\[\mathbb E[D_n^2] = \mathbb{E} \left[\sum_{k=1}^{n-1}(2k-1) \mathbf 1_{\{D_n\ge k\}} \right] = \sum_{k=1}^{n-1}(2k-1) \mathbb{E}[\mathbf{1}_{\{D_n\ge k\}}].\]
Recall that, for an event $X$, we have $\displaystyle \mathbb{E}[\mathbf{1}_X] = \operatorname{Pr}(X)$. Hence,
\[\mathbb E[D_n^2] = \sum_{k=1}^{n-1}(2k-1)\Pr(D_n\ge k).\]

\noindent Using the previously established bound in (\ref{bound for Pr(D_n>=k)}):
\[\Pr(D_n\ge k)\le\frac{\rho^k}{1-\rho},\] we get
\begin{align*}
\mathbb E[D_n^2]
&\leq \frac{1}{1-\rho}
\sum_{k=1}^{n-1}(2k-1)\rho^k\\
&\leq \frac{1}{1-\rho}
\sum_{k=1}^{\infty}(2k-1)\rho^k\\
&= \frac{\rho(1+\rho)}{(1-\rho)^3}.
\end{align*}

\noindent Thus, from (\ref{bound estimate for variance}), we have
\begin{align}
     V_{\mu_{\mathbf p}}(\mathcal C_n) \leq \frac{\rho}{(1-\rho)^2}.
\end{align}

\noindent Independence of the Bernoulli coordinates gives us
\[H_{\mu_{\mathbf p}}(\mathcal C_n) = nH(\mathbf p) \qquad \text{and} \qquad h_{\mu_{\mathbf p}}(\sigma) = H(\mathbf p)>0.\]

\noindent From (\ref{A(C_n) = n-E[D(n)]}), $\displaystyle A_{\mu_{\mathbf p}}(\mathcal C_n) = n-\mathbb E[D_n]$. Because $\mathbb E[D_n]$ is uniformly bounded, we then have
\begin{align}\label{limit of A(C_n)/(1+H(C_n)) = 1/H(p)}
\lim_{n\to\infty} \frac{A_{\mu_{\mathbf p}}(\mathcal C_n)}
{1+H_{\mu_{\mathbf p}}(\mathcal C_n)} = \lim_{n\to\infty}
\frac{n-\mathbb E[D_n]} {1+nH(\mathbf p)} =\frac{1}{H(\mathbf p)}.
\end{align}

\noindent Similarly, since the variance is uniformly bounded
and $H(\mathbf p)>0$, we have
\begin{align}\label{limit variance = 0}
    \lim_{n\to\infty} \frac{V_{\mu_{\mathbf p}}(\mathcal C_n)}
{(1+H_{\mu_{\mathbf p}}(\mathcal C_n))^2} = 0.
\end{align}

\noindent Finally, consider the canonical partition
$\displaystyle \mathcal{C} = \{[0], [1], \ldots, [q-1]\}$ of the one-sided Bernoulli shift. Since the alphabet is finite and the cylinder sets generate the product $\sigma$-algebra,
$\mathcal{C}$ is a finite forward generator. Moreover, by Lemma \ref{lemma: Bernoulli Return}, we have $\displaystyle M_{\mu_{\mathbf p}}(\mathcal C_n)=n$. Thus,
\[M_{\mu_{\mathbf p}}(\mathcal C_n)=O(n)
\quad\text{as }n\to\infty.\]
(Recall that this means there exist $C>0$ and $N\in \mathbb{N}$ such that $\displaystyle M_{\mu_{\mathbf p}}(\mathcal C_n) <nC$ for all $n>N$. In other words, the ratio $\displaystyle \frac{M_{\mu_{\mathbf p}}(\mathcal C_n)}{n}$ is bounded for all sufficiently large $n$.)

\noindent So, $\mathcal{C}\in \mathscr A_{\mu_{\textbf{p}}}(\sigma)$. That is, $\mathcal{C}$ is linearly return-time admissible. Since the global invariants are defined as suprema over all such partitions, the values obtained for $\mathcal{C}$, in (\ref{limit of A(C_n)/(1+H(C_n)) = 1/H(p)}) and (\ref{limit variance = 0}), give us the desired global lower bounds:
\[\Lambda_{\mu_{\mathbf p}}(\sigma)\ge\frac1{H(\mathbf p)} \qquad \text{and} \qquad \Lambda_{V,\mu_{\mathbf p}}(\sigma)\ge0.\]
This completes the proof of Theorem \ref{theorem: Canonical Bernoulli calculations}.
\end{proof}

\begin{example}[Uniform Bernoulli shifts]
\label{ex:bernoulli-uniform}
Consider the one-sided Bernoulli shift $(\Sigma_q^+,\mu_{\mathbf p},\sigma)$ on the alphabet $\mathcal A=\{0,1,\ldots,q-1\}$, where $q\ge 2$, with the uniform probability distribution vector
$\displaystyle \mathbf p=\left(\frac1q,\ldots,\frac1q\right)$.
Thus, each symbol occurs independently with probability $1/q$.

\noindent The Shannon entropy of this distribution is
\[H(\mathbf p) = -\sum_{i=1}^{q}p_i\ln p_i =-\sum_{i=1}^{q}\frac1q\ln\frac1q =-q\cdot\frac1q\ln\frac1q = \ln q.\]

\noindent Let $\displaystyle \mathcal{C} = \{[0], [1], \ldots, [q-1]\}$ be the canonical partition into one-symbol cylinders, and let
\[\mathcal C_n = \bigvee_{j=0}^{n-1}\sigma^{-j}\mathcal C.\]
Since every word of length $n$ has probability $q^{-n}$, the entropy of $\mathcal C_n$ is
\[ H_{\mu_{\mathbf p}}(\mathcal C_n) = -\sum_{w\in\mathcal A^n}
\frac{1}{q^n}\ln(q^{-n}) = -q^n\cdot q^{-n}\ln(q^{-n}) = n\ln q.\]

\noindent From (\ref{Average = Expectation}) and (\ref{A(C_n) = n-E[D(n)]}), the expected smallest period of a random word of length $n$ satisfies
\begin{align}\label{E[per (W_n)] = n - E[D_n]}
    \mathbb{E}[\operatorname{per}(W_n)] = A_{\mu_{\mathbf p}}(\mathcal C_n) = n-\mathbb E[D_n],
\end{align}
where $\mathbb E[D_n]$ is bounded independently of $n$. Thus,
\[L_{\mu_{\mathbf p}}(\sigma,\mathcal C) :=\lim_{n\to\infty}
\frac{A_{\mu_{\mathbf p}}(\mathcal C_n)}{1+H_{\mu_{\mathbf p}}(\mathcal C_n)} = \lim_{n\to\infty}\frac{n-\mathbb E[D_n]}{1+n\ln q} = \frac1{\ln q}.\]

\noindent Similarly, since the variance $V_{\mu_{\mathbf p}}(\mathcal C_n)$ is uniformly bounded, we obtain
\[L_{V,\mu_{\mathbf p}}(\sigma,\mathcal C) = \lim_{n\to\infty}
\frac{V_{\mu_{\mathbf p}}(\mathcal C_n)}{(1+H_{\mu_{\mathbf p}}(\mathcal C_n))^2} = 0.\]

\noindent Since $\mathcal C$ is linearly return-time admissible, these calculations also give
\[\Lambda_{\mu_{\mathbf p}}(\sigma)
\ge \frac1{\ln q} \qquad \text{and} \qquad \Lambda_{V,\mu_{\mathbf p}}(\sigma)\ge 0.\]
However, these inequalities are not necessarily equalities, since the global invariants are defined by taking suprema over all linearly return-time admissible partitions.
\end{example}

\begin{example}[A Nonuniform Bernoulli Measure]
\label{ex:bernoulli-biased}
Consider the one-sided binary Bernoulli shift
$(\Sigma_2^+,\mu_p,\sigma)$ with alphabet
$\mathcal A=\{0,1\}$ and probability distribution vector
$\displaystyle \textbf{p}=(p,1-p)$, where $0<p<1$. Let $\mathcal C = \{[0], [1]\}$ be the canonical partition into one-symbol cylinders. Thus, the symbols occur independently, with probabilities
\[\mu_p([0])=p \qquad \text{and} \qquad \mu_p([1])=1-p.\]
Let
\[\mathcal C_n =\bigvee_{j=0}^{n-1}\sigma^{-j}\mathcal C.\]
The Shannon entropy of the probability vector \textbf{p} is
\[H(p) = -p\ln p - (1-p)\ln(1-p).\]

\noindent Since the Bernoulli coordinates are independent,
the entropy of the partition $\mathcal C_n$ is
\[H_{\mu_p}(\mathcal C_n)=nH(p) = -n\left(p\ln p + (1-p)\ln(1-p) \right).\]

\noindent From (\ref{E[per (W_n)] = n - E[D_n]}), we have
\[A_{\mu_p}(\mathcal C_n)
=\mathbb E[\operatorname{per}(W_n)]
=n-\mathbb E[D_n].\]
Here, $\displaystyle D_n=n-\operatorname{per}(W_n)$ is the length of the longest proper border of the random word $W_n$. Since $\mathbb E[D_n]$ is bounded independently
of $n$, we obtain
\[\displaystyle L_{\mu_p}(\sigma,\mathcal C) = \lim_{n\to\infty}\frac{A_{\mu_p}(\mathcal C_n)}{1+H_{\mu_p}(\mathcal C_n)} = \lim_{n\to\infty}
\frac{n-\mathbb E[D_n]}{1+nH(p)} = \frac1{H(p)} = \frac{-1}{p\ln p+(1-p)\ln(1-p)}.\]

\noindent Similarly, Theorem \ref{theorem: Canonical Bernoulli calculations} shows that $V_{\mu_p}(\mathcal C_n)$ is bounded
independently of $n$. Consequently,
\[ L_{V,\mu_p}(\sigma,\mathcal C) = \lim_{n\to\infty}
\frac{V_{\mu_p}(\mathcal C_n)}{(1+H_{\mu_p}(\mathcal C_n))^2} =0.\]

\noindent To compare different values of $p$, notice that
\[H'(p)=\ln\left(\frac{1-p}{p}\right) \qquad \text{and} \qquad H''(p)=-\frac1p-\frac1{1-p}<0.\]

\noindent Thus, $H(p)$ is strictly concave and attains
its unique maximum at $p=1/2$. Thus,
\[L_{\mu_p}(\sigma,\mathcal C) = \frac{1}{H(p)}\ge \frac{1}{H\left( \frac{1}{2} \right)} =\frac{1}{\ln 2},\]
with equality if and only if $p=1/2$.

\noindent Therefore, among binary Bernoulli measures, the fair (non-biased) measure minimizes the normalized mean-return invariant associated with the canonical partition. This comparison illustrates the entropy normalization and does not imply that biased Bernoulli measures have faster actual return times.

\noindent Finally, since $\mathcal C$ is a finite forward generator and $\displaystyle M_{\mu_p}(\mathcal C_n)=n$, it is linearly return-time admissible: $\mathcal{C}\in \mathscr A_{\mu_{\textbf{p}}}(\sigma)$. Hence, the definitions of the global invariants give
\[\Lambda_{\mu_p}(\sigma) \ge\frac1{H(p)} \qquad \text{and} \qquad \Lambda_{V,\mu_p}(\sigma)\ge0.\]
Notice that these are lower bounds, not necessarily equalities.
\end{example}

\begin{example}[Return times within one block partition]
\label{ex:bernoulli-words}
For the fair binary shift and the partition into words of
length three: 
\[\mathcal{C}_3 := \{[000], [001], [010], [011], [100], [101], [110], [111]\}.\] 
The eight words (the eight atoms of the partition) have minimum overlap times:
\[ \begin{array}{|c|c|c|c|c|c|c|c|c|}
 \hline w&000&001&010&011&100&101&110&111\\ 
 \hline
 \tau_\mu([w])&1&3&2&3&3&2&3&1 \\
 \hline
 \end{array}
\]
Consequently,
\begin{align*}
    A_\mu(\mathcal C_3) &= \frac{1}{8}(1+3+2+3+3+2+3+1) = \frac{18}{8}=\frac94, \quad \text{and}\\
     V_\mu(\mathcal C_3) &=\frac{1}{8}\left(2\left(1-\frac94\right)^2 + 2\left(2-\frac94\right)^2 + 4\left(3-\frac94\right)^2\right) = \frac{11}{16}.
\end{align*}

\noindent This illustrates why the variance is nonzero at a finite observation length even though its entropy-normalized limit is zero.
\end{example}

\begin{proposition}[Finite-Window Partitions]
\label{prop:finite-window-return}
Let $(\Sigma,\mu_{\mathbf p},\sigma)$ be a one-sided Bernoulli shift over a finite alphabet with probability distribution vector $\mathbf p$. Suppose that $\mathcal {R}$ is a finite measurable
partition determined by coordinates $0,\ldots,r-1$, where $r\ge 1$. Then, for every $n\ge 1$, we have
\[ M_{\mu_{\mathbf p}}(\mathcal R_\sigma^{(n)})
\le n+r-1. \]

\noindent For a two-sided Bernoulli shift, if $\mathcal R$ is determined by coordinates $-r,\ldots,r$, then
\[M_{\mu_{\mathbf p}}(\mathcal R_\sigma^{(n)})
\le n+2r.\]
\end{proposition}

\begin{proof}
First, we consider the one-sided case. We say $\mathcal {R}$ is determined by coordinates $0,\ldots,r-1$ if any two sequences that agree on these coordinates belong to the same atom of $\mathcal R$. In other words, knowing the first $r$ symbols of a sequence is sufficient to determine which atom of
$\mathcal R$ contains it. Recall that $\displaystyle \mathcal R_\sigma^{(n)} = \bigvee_{j=0}^{n-1}\sigma^{-j}( \mathcal{R})$. Since $\sigma^{-j}(\mathcal{R})$ depends only on coordinates $j,\ldots,j+r-1$, every atom of $\mathcal R_\sigma^{(n)}$ is
determined by coordinates $\displaystyle 0,\ldots,n+r-2$.

Let $B$ be an atom of $\mathcal R_\sigma^{(n)}$ with
$\mu_{\mathbf p}(B)>0$, and let $\displaystyle L=n+r-1$. Then, $B$ depends only on coordinates $0,\ldots,L-1$, while
$\sigma^{-L}(B)$ depends only on coordinates $L,\ldots,2L-1$.

\noindent These coordinate sets are disjoint. Since the coordinates of a Bernoulli shift are independent, we have
\[\mu_{\mathbf{p}}(B\cap\sigma^{-L}(B)) = \mu_{\mathbf{p}}(B)\mu_{\mathbf{p}}(\sigma^{-L}(B)).\]
Because $\mu_{\mathbf{p}}$ is $\sigma$-invariant, we have $\displaystyle \mu_{\mathbf{p}}(\sigma^{-L}(B)) = \mu_{\mathbf{p}}(B)$. Therefore,
\[\mu_{\mathbf p}(B\cap\sigma^{-L}B)
=\mu_{\mathbf p}(B)^2>0.\]

\noindent By the definition of the return time, $\displaystyle \tau_{\mu_{\mathbf p}}(B)\le L=n+r-1$. Taking the maximum over all
positive-measure atoms gives
\[M_{\mu_{\mathbf p}}(\mathcal R_\sigma^{(n)})
\le n+r-1.\]

\noindent For the two-sided case, the partition $\displaystyle \mathcal{R}_\sigma^{(n)}$ is determined by coordinates $\displaystyle -r,\ldots,n+r-1$. This window has length $n+2r$.
Shifting it by $n+2r$ times gives us a disjoint coordinate window. The same independence argument, therefore, gives
\[M_{\mu_{\mathbf p}}(\mathcal R_\sigma^{(n)})
\le n+2r.\]
\end{proof}

\begin{remark}
The preceding proposition applies to partitions whose atoms are determined by a finite number of coordinates. For such partitions, shifting an atom by a sufficiently large number of coordinates produces two events depending on disjoint coordinate sets. The independence of Bernoulli coordinates then guaranties a positive-measure intersection, giving us a linear bound on the maximal return time.

\noindent However, an arbitrary finite measurable generator need not be determined by any fixed finite coordinate window. Its atoms
may depend on infinitely many coordinates, so shifting an atom by a finite number of positions does not necessarily produce
events depending on disjoint coordinate sets, causing the
 independence argument used above to no longer be directly applicable.

\noindent Although measurable sets can be approximated in measure by finite-window sets, such approximations do not automatically provide a uniform linear bound on the maximal return times of the successive refinements. Therefore, the preceding proposition does
not establish that every finite measurable generator of a Bernoulli shift is linearly return-time admissible.
\end{remark}

\subsection{Finite atomic systems}
\label{subsection:finite-atomic-measure}

Let $X$ be finite, let $\mu(\{x\})>0$ for every $x\in X$,
and let $f:X\to X$ preserve $\mu$. Here $\mu$ is a probability measure. Then, $f$ is a permutation. Let its disjoint periodic be denoted by $O_1,\ldots,O_r$. Let $m_i$ and $a_i = \mu(O_i)$ be the lengths and the masses of $O_i$, respectively, for $i=1, \ldots, r$. So, because $f$ is measure preserving, every point of $O_i$ has mass $a_i/m_i$.

\begin{theorem}[Exact global values for finite permutations]
\label{Theorem: Exact global values for finite permutations}
For a finite full-support probability-preserving permutation, we have
\[ \Lambda_\mu(f)= \frac{\displaystyle \sum_{i=1}^r a_i m_i}
{\displaystyle 1-\sum_{i=1}^r a_i\log\left( \frac{a_i}{m_i} \right)} \qquad \text{and} \qquad 
\Lambda_{V,\mu}(f)=
\frac{\displaystyle\sum_{i=1}^r a_i
\left(m_i-\sum_{j=1}^r a_jm_j\right)^2}
{\displaystyle \left(1-\sum_{i=1}^r a_i\log\left(\frac{a_i}{m_i} \right)\right)^2}.\]
\end{theorem}
\begin{proof}
Let $\mathcal S=\{\{x\}:x\in X\}$ be the singleton partition of $X$. Clearly, it is generating and admissible, because
$\mathcal S_f^{(n)}=\mathcal S$ and $\displaystyle M_\mu(\mathcal S)=\max_{1\leq i \leq |X|} m_i<\infty$.

\noindent Now let $\displaystyle \mathcal P\in\mathscr{A}_\mu(f)$ be an arbitrary admissible generator. Notice that, since $X$ is finite, the sequence of refinements $\displaystyle \left(\mathcal{P}_f^{(n)} \right)_{n=1}^\infty$ eventually stabilizes. That is, there exists $N\geq 1$ such that $\displaystyle \mathcal {P}_f^{(n)}=\mathcal {P}_f^{(N)} = \mathcal{S}$ for all $n\geq N$. (Since $\mathcal {P}$ is forward generating, this limiting partition must be the singleton partition $\mathcal S$.) Hence, we get
\[L_\mu(f,\mathcal P)
=\frac{A_\mu(\mathcal {S})}{1+H_\mu(\mathcal {S})} \quad \text{and} \quad L_{V,\mu}(f,\mathcal {P})
=\frac{V_\mu(\mathcal {S})}{(1+H_\mu(\mathcal {S}))^2}.\]

\noindent It is easy to see that each $x\in O_i$ has $\tau_\mu(\{x\})=m_i$ and $\displaystyle \mu(\{x\})=\frac{a_i}{m_i}$. Therefore, we have
\begin{align*}
A_\mu(\mathcal S)&=\sum_i a_i m_i,\\
H_\mu(\mathcal S)&=-\sum_i a_i\log\left( \frac{a_i}{m_i} \right), \quad \text{and}\\
V_\mu(\mathcal {S})&=\sum_i a_i
\left(m_i-\sum_j a_jm_j\right)^2.
\end{align*}
Because $\displaystyle \mathcal{P} \in \mathscr {A}_\mu (f)$ is arbitrary, these values are independent of the choice of admissible generator. Since the singleton partition $\mathcal S$ is itself admissible, taking the supremum over all admissible
generators gives
\[\Lambda_\mu(f) = \frac{A_\mu(\mathcal S)} {1+H_\mu(\mathcal S)} = \frac{\displaystyle \sum_{i=1}^r a_i m_i}
{\displaystyle 1-\sum_{i=1}^r a_i\log\left( \frac{a_i}{m_i} \right)}\]
and
\[\Lambda_{V,\mu}(f) =\frac{V_\mu(\mathcal S)} {(1+H_\mu(\mathcal S))^2} = \frac{\displaystyle\sum_{i=1}^r a_i
\left(m_i-\sum_{j=1}^r a_jm_j\right)^2}
{\displaystyle \left(1-\sum_{i=1}^r a_i\log\left(\frac{a_i}{m_i} \right)\right)^2}.\]
\end{proof}

\noindent The following result follows immediately.
\begin{corollary}[One periodic orbit case]
\label{Corollary: One periodic orbit case}
If $f$ cyclically permutes $m$ points with uniform measure. Then,
\[\Lambda_\mu(f)=\frac{m}{1+\log m} \qquad 
\text{and} \qquad \Lambda_{V,\mu}(f)=0.\]
\end{corollary}

\begin{example}[Equal entropy but distinct return-time invariants]
\label{example: same entropy but different return times}
Let $C_2$ and $C_3$ be the cyclic permutations on two and
three points, respectively, each equipped with a uniform measure. Both have Kolmogorov-Sinai entropy zero. However,
\[\Lambda_{\mu_2}(C_2)=\frac2{1+\log2} \ne\frac3{1+\log3} = \Lambda_{\mu_3}(C_3).\]
Thus, $\Lambda_\mu$ is not determined by the Kolmogorov-Sinai entropy, even within finite systems.
\end{example}

\begin{example}[Positive global variance]
\label{example: Positive global variance}
Let $X$ consist of a fixed point and a two-cycle, with uniform measure on the three points. Then, we have
$A_\mu(\mathcal S)=5/3$, $V_\mu(\mathcal S)=2/9$, and $H_\mu(\mathcal S)=\log3$. Thus,
\[\Lambda_\mu(f)=\frac{5}{3(1+\log3)} \qquad \text{and} \qquad \Lambda_{V,\mu}(f)=\frac{2}{9(1+\log3)^2}>0.\]
In particular, the global variance invariant need to vanish.
\end{example}

\begin{remark}
All finite permutations have zero Kolmogorov-Sinai
entropy, regardless of their cycle lengths or
cycle-mass distributions. However, the formulas in
Theorem~\ref{Theorem: Exact global values for finite permutations}
show that the mean invariant retains information
about the probability-weighted average period,
while the variance invariant reflects the
probability-weighted variation in period lengths.
Thus, these invariants can distinguish finite
systems that have the same zero entropy.
Nevertheless, they do not necessarily provide
a complete classification of finite permutations.
\end{remark}

\begin{example}[Contrasting mean and variance behaviors]
\label{example: Contrasting mean and variance behaviors}
Let $f$ be the cyclic permutation of two points with uniform
measure. Let $g$ be the permutation of three points consisting
of a fixed point and a two-cycle, with uniform measure. Then,
\[ \Lambda_\mu(f)=\frac2{1+\log2},  \qquad \Lambda_{V,\mu}(f)=0,
\]
whereas
\[ \Lambda_\nu(g)=\frac{5}{3(1+\log3)},  \qquad  \Lambda_{V,\nu}(g)=\frac{2}{9(1+\log3)^2}>0.\]

\noindent The values of the mean invariant in this particular pair
are not equal. This example shows that nonzero variance can
occur in finite systems even though it vanishes in a single
cycle. Finding a pair with an equal global mean but different
global variance is a separate problem.
\end{example}

\subsection{Approximation by finite periodic-orbit systems}
\label{subsection:periodic-approximation}

Let $q\in \mathbb{N}$, with $q\geq 2$, be fixed, and consider the full shift on $q$ symbols. For each $N\in \mathbb{N}$, let
\begin{align}\label{Finite Systems}
    X_N:=\bigcup_{k=1}^N\operatorname{Fix}(\sigma^k) \qquad \text{and} \qquad f_N:=\sigma|_{X_N}.
\end{align}

\noindent We quip $X_N$ with the uniform probability measure $\mu_N$.
Every one-sided periodic sequence has a unique two-sided periodic extension, and this correspondence preserves the shift action.
Consequently, the restrictions of the one-sided and two-sided full shifts to their periodic points of least period at most $N$ are isomorphic as finite dynamical systems. In other words, whether we start with the one-sided full shift or the two-sided full shift, the resulting finite dynamical systems $\displaystyle (X_N, f_N)$ are equivalent. \\

\noindent Let $P_q(d)$ be the number of points of least period exactly $d$. Let
\[s_N:=|X_N|=\sum_{d=1}^N P_q(d) \qquad \text{and} \qquad \widetilde{d}_N:=\frac{\displaystyle\sum_{d=1}^N dP_q(d)}{s_N}.\]

\noindent Recall that the number-theoretic M\"obius function
$\operatorname{mob}:\mathbb N\to\{-1,0,1\}$
is defined by
\[\operatorname{mob}(n)=
\begin{cases}
1, & n=1,\\
(-1)^r, & \text{if }n\text{ is a product of }
r\text{ distinct primes},\\
0, & \text{if }n\text{ is divisible by the
square of a prime}.
\end{cases}\]

\noindent We will need the following useful formula to calculate approximations by periodic orbits (Theorem \ref{Theorem:Periodic Orbit Approximation}).
\begin{theorem}[M\"obius Inversion Formula]
Let $F,G:\mathbb N\to\mathbb C$ be arithmetic
functions. Then,
\[F(n)=\sum_{d\mid n}G(d)\]
for every $n\in\mathbb {N}$ if and only if
\[G(n)=\sum_{d\mid n}\operatorname{mob}(d)F\left(\frac{n}{d}\right),\]
where $\operatorname{mob}:\mathbb{N}\to\{-1,0,1\}$ is the number-theoretic M\"obius function defined above.
\end{theorem}

\noindent For the sake of completeness, we will provide a short proof here. For more details, one may refer to \cite{Apostol1976} (Sections 2.2 and 2.7), \cite{NISTDLMF} (Sections 27.2 and 27.5)

\begin{proof}
From the definition, the M\"obius function satisfies
\[\sum_{d\mid n}\operatorname{mob}(d) = 
\begin{cases}
1, & n=1,\\
0, & n>1.
\end{cases}\]
The reason for the case $n>$ is:
\[\sum_{d\mid n}\operatorname{mob}(d) = \sum_{k=0}^r \binom{r}{k}(-1)^k = (1-1)^r=0.\]

\noindent ($\Rightarrow$) Suppose that $\displaystyle F(n)=\sum_{d\mid n}G(d)$. Then,

\begin{align*}
\sum_{d\mid n}\operatorname{mob}(d)
F\left(\frac nd\right)
&= \sum_{d\mid n}\operatorname{mob}(d) \sum_{k\mid (n/d)}G(k)\\
&= \sum_{k\mid n}G(k) \sum_{d\mid (n/k)}\operatorname{mob}(d)\\
&= \sum_{k\neq n}G(k) \cdot 0 + G(n) \cdot 1\\
&=G(n).
\end{align*}

\noindent ($\Leftarrow$) By using the same arguments and calculations, one can obtain the converse.
\end{proof}

\noindent Let $P_q(d)$ denote the number of periodic points of the full shift on $q$ symbols whose least period is exactly $d$. Because every point fixed by $\sigma^n$ has least period dividing $n$, we have
\[q^n=\sum_{d\mid n}P_q(d).\]

\noindent The M\"obius inversion formula gives us
\begin{align}\label{Mobius Inversion for shift space}
    P_q(d)=\sum_{k\mid d}\operatorname{mob}(k)q^{d/k}.
\end{align}

\noindent Since every periodic orbit of length $d$ contains
exactly $d$ distinct points, it follows that the number of distinct periodic orbits of length $d$ is:
\begin{align}\label{the number of distinct periodic orbits}
    \frac{P_q(d)}{d} = \frac{1}{d}\sum_{k\mid d} \operatorname{mob}(k)q^{d/k}.
\end{align}

\begin{theorem}[Periodic Orbit Approximation]
\label{Theorem:Periodic Orbit Approximation}
For the finite systems defined in \eqref{Finite Systems}, we have
\begin{align*}
    \lim_{N\to\infty}\Lambda_{\mu_N}(f_N) &= \lim_{N\to\infty} \frac{\widetilde d_N}{1+\log s_N} = \frac{1}{\log q}, \quad \text{and}\\
    \\
    \lim_{N\to\infty}\Lambda_{V,\mu_N}(f_N) &= \lim_{N\to\infty}
\frac{\displaystyle \frac{1}{s_N}\sum_{d=1}^{N} P_q(d)(d-\widetilde d_N)^2 }{ (1+\log s_N)^2} = 0.
\end{align*}
\end{theorem}

\begin{proof}
First, we derive the exact formulas for the invariants $(X_N,f_N,\mu_N)$.

\noindent Notice that, because $X_N$ consists of periodic points,
$f_N$ is a permutation of the finite set $X_N$. Moreover, $\mu_N$ is the uniform probability measure. So,
\[\mu_N(\{x\})=\frac{1}{s_N} \qquad\text{for every }x\in X_N. \]

\noindent Let $\mathcal{S}_N$ denote the singleton partition of $X_N$. Since $\mathcal{S}_N$ is a generating partition, Theorem
\ref{Theorem: Exact global values for finite permutations} gives us
\begin{align}\label{Lambda_mu and Lambda_V}
    \Lambda_{\mu_N}(f_N) = \frac{A_{\mu_N}(\mathcal S_N)}{1+H_{\mu_N}(\mathcal S_N)} \qquad \text{and} \qquad \Lambda_{V,\mu_N}(f_N) = \frac{V_{\mu_N}(\mathcal S_N)}{(1+H_{\mu_N}(\mathcal S_N))^2}.
\end{align}

\noindent By definition, $P_q(d)$ is the number of points of least period exactly $d$. Thus, among the $s_N$ number of points of $X_N$, exactly $P_q(d)$ have least period $d$. Since every point has measure $\displaystyle \frac{1}{s_N}$, the total measure of these points of least period $d$ is $\displaystyle \frac{P_q(d)}{s_N}$. Notice that for each $x\in X_N$, the first return time of the singleton $\{x\}$ is equal to the least period of $x$. Therefore, by grouping the points according to their
least periods, we get
\begin{align}\label{A_mu_N(S_N) = d_N}
    A_{\mu_N}(\mathcal{S}_N) = \sum_{x\in X_N}\mu_N(\{x\})\operatorname{per}(x) = \frac{1}{s_N}
\sum_{d=1}^{N}dP_q(d) = \widetilde{d}_N.
\end{align}

\noindent Similarly, for the variance, we get
\begin{align} \label{V_mu_N(S_N) = sum}
    V_{\mu_N}(\mathcal{S}_N) = \sum_{x\in X_N}\mu_N(\{x\}) \bigl(\operatorname{per}(x) - \widetilde{d}_N\bigr)^2 = \frac{1}{s_N}\sum_{d=1}^{N} P_q(d)(d-\widetilde{d}_N)^2.
\end{align}

\noindent Furthermore, since $\mathcal{S}_N$ consists of $s_N$ number of atoms, each of measure $1/s_N$, its (measure theoretic) entropy is
\begin{align}\label{Entropy of S_N}
    H_{\mu_N}(\mathcal{S}_N) = -\sum_{x\in X_N} \frac{1}{s_N}\log\left(\frac{1}{s_N}\right) = \log s_N.
\end{align}

\noindent Thus, by (\ref{Lambda_mu and Lambda_V}), (\ref{A_mu_N(S_N) = d_N}), (\ref{V_mu_N(S_N) = sum}), and (\ref{Entropy of S_N}), we get
\begin{align}\label{Lambda(f_N)}
    \Lambda_{\mu_N}(f_N) = \frac{\widetilde d_N}{1+\log s_N}
\end{align}
and
\begin{align}\label{Lambda_V}
    \Lambda_{V,\mu_N}(f_N) = \frac{\displaystyle \frac{1}{s_N}\sum_{d=1}^{N}P_q(d)(d-\widetilde d_N)^2}{(1+\log s_N)^2}.
\end{align}

\noindent We now establish the asymptotic behavior of these expressions as $N\to\infty$. By the M\"obius inversion formula, we got \eqref{Mobius Inversion for shift space}: $\displaystyle P_q(d) = \sum_{k\mid d}\operatorname{mob}(k)q^{d/k}$. By separating the term corresponding to $k=1$, we get
\[P_q(d) = q^d + \sum_{\substack{k\mid d\\k>1}}\operatorname{mob}(k)q^{d/k}.\]

\noindent Because $k\geq 2$ in the remaining sum, we have $d/k\leq d/2$. Also, from definition, $|\operatorname{mob}(k)|\leq 1$, and there are at most $d$ number of divisors of $d$. Hence, $\displaystyle |P_q(d)-q^d| \leq dq^{d/2}$, and therefore,
\[\lim_{d\rightarrow \infty }\frac{P_q(d)}{q^d} = 1.\]
Thus, there exists a constant $c_q>0$ such that $\displaystyle P_q(N)\geq c_q q^N$, for all $N\geq 1$. On the other hand, since $P_q(d)\leq q^d$, we have
\[s_N=\sum_{d=1}^{N}P_q(d) \leq\sum_{d=1}^{N}q^d \leq\frac{q}{q-1}q^N.\]

\noindent Because $\displaystyle s_N\geq P_q(N)$, we then get
\[c_q q^N\leq s_N \leq\frac{q}{q-1}q^N, \qquad \text{for all } N\geq 1.\]

\noindent This implies $\displaystyle s_N\asymp q^N$. By taking logarithms, we then obtain
\begin{align}\label{log S_N}
    \log s_N=N\log q+O(1).
\end{align}

\noindent Now, we show that the average period $\widetilde d_N$ differs from $N$ by a uniformly bounded amount. Consider the period deficit $N-d$. Its expected value under the weights
$P_q(d)/s_N$ is:
\[N-\widetilde d_N = N-\frac{\displaystyle \sum_{d=1}^N dP_q(d)}{s_N} =  \frac{1}{s_N} \sum_{d=1}^{N}(N-d)P_q(d).\]
Because $N-d\geq 0$, this implies $N-\widetilde{d}_N \geq 0$. Let $j=N-d$. Then,
\[N-\widetilde d_N = \sum_{j=0}^{N-1} j\frac{P_q(N-j)}{s_N}. \]

\noindent Notice that, because $P_q(N-j)\leq q^{N-j}$ and $s_N\geq c_q q^N$, we get $\displaystyle \frac{P_q(N-j)}{s_N}
\leq \frac{1}{c_q}q^{-j}$. Thus,
\[0\leq N-\widetilde d_N \leq \frac{1}{c_q} \sum_{j=0}^{N-1}jq^{-j} \leq \frac{1}{c_q} \sum_{j=0}^{\infty}jq^{-j} <\infty.\]

\noindent (The series converges because $q\geq 2$.) Hence, $\displaystyle \widetilde d_N=N+O(1)$. This together with (\ref{log S_N}) gives us
\[\lim_{N\rightarrow \infty}\Lambda_{\mu_N}(f_N) = \lim_{N\rightarrow \infty} \frac{\widetilde d_N}{1+\log s_N} = \lim_{N\rightarrow \infty} \frac{N+O(1)} {N\log q+O(1)} = \frac{1}{\log q}.\]

\noindent It remains to consider the variance invariant. Notice that
\[\mathbb{E}[(N-d)^2]=\frac{1}{s_N} \sum_{d=1}^{N}(N-d)^2P_q(d) = \sum_{j=0}^{N-1}
j^2\frac{P_q(N-j)}{s_N} \leq \frac{1}{c_q} \sum_{j=0}^{\infty}j^2q^{-j}<\infty.\]
This says that the second moment of the period deficit is uniformly bounded. Recall the mean minimizes the mean square deviation. In fact, since $\widetilde {d}_N$ is the weighted mean
of the period lengths, we get
\[ \frac{1}{s_N}\sum_{d=1}^{N} P_q(d)(d-N)^2 = \frac{1}{s_N}\sum_{d=1}^{N} P_q(d)(d-\widetilde d_N)^2 + (N-\widetilde {d}_N)^2.\]
This implies
\[\frac{1}{s_N}\sum_{d=1}^{N} P_q(d)(d-\widetilde d_N)^2 \leq \frac{1}{s_N}\sum_{d=1}^{N} P_q(d)(d-N)^2 = O(1).\]

\noindent Therefore,
\[0 \leq \lim_{N \rightarrow \infty} \Lambda_{V,\mu_N}(f_N) =\lim_{N \rightarrow \infty} \frac{\displaystyle \frac{1}{s_N}\sum_{d=1}^{N}P_q(d)(d-\widetilde d_N)^2}{(1+\log s_N)^2} = \lim_{N \rightarrow \infty}\frac{O(1)} {(N\log q+O(1))^2} = 0.\]

\noindent This completes the proof of the theorem.
\end{proof}

\begin{remark}
The limiting value $\displaystyle \frac{1}{\log q}$ coincides with
the value of the mean invariant associated with the canonical partition of the uniform Bernoulli $q$-shift. That is, if we let $\displaystyle \textbf{p}= \left(\frac{1}{q}, \ldots, \frac{1}{q} \right)$ in Theorem \ref{theorem: Canonical Bernoulli calculations}, then we get $\displaystyle L_{\mu_\textbf{p}} (\sigma, \mathcal{C})=\frac{1}{\log q}$. However, this convergence concerns a sequence of finite measure-preserving
systems, each having zero Kolmogorov-Sinai entropy. In particular, it does not establish that the global invariant of the infinite Bernoulli shift equals $\displaystyle \frac{1}{\log q}$.
\end{remark}

\section{Further Questions}
\label{Section: Further Questions}

The preceding results raise questions concerning the
necessity of admissibility, exact global values, and the
relationship between topological and measure-theoretic
return-time statistics.

\begin{question}[One-sided full shifts]
\label{question:topological one sided}
For the one-sided full shift $\sigma_q$ on $q\ge2$ symbols, the canonical cylinder generator gives $\displaystyle L_A(\sigma_q,\mathcal C_1)=\frac{1}{\log q}$. (See (\ref{L_A(sigma_q, C_1)})). Does the global topological invariant also satisfy it? That is, do we have
\[\Lambda_A(\sigma_q)=\frac1{\log q}?\]
\end{question}

\begin{question}[Topological variance for full shifts]
\label{question:topological variance}
The canonical cylinder generator satisfies
$L_V(\sigma_q,\mathcal C_1)=0$.
Is $\Lambda_V(\sigma_q)=0$, or can a different
admissible generating cover yield positive normalized variance?
\end{question}

\begin{question}[Two-sided topological full shifts]
\label{question:topological two sided}
For the two-sided full shift, using the symmetric-join
convention adopted in the topological part of the paper,
the canonical generator has
\[\widehat{L}_A(\sigma_q,\mathcal C_0)=\frac1{\log q} \qquad \text{and} \qquad \qquad \widehat{L}_V(\sigma_q,\mathcal C_0)=0.\]
(See Theorem \ref{Theorem: L_A(sigma,C_0)=1/log q} and Proposition \ref{proposition variance symmetric finite}.) Are the corresponding global quantities equal to these canonical values? That is, do we have
\[\widehat{L}_A(\sigma_q)=\frac1{\log q} \qquad \text{and} \qquad \qquad \widehat{L}_V(\sigma_q)=0\ ?\]
\end{question}

\begin{question}[Necessity of measure-theoretic admissibility]
\label{question:measure admissibility}
Does there exist a finite measurable generating partition
$\mathcal{P}$ of a nontrivial Bernoulli shift such that
\[\limsup_{n\to\infty} \frac{M_\mu(\mathcal{P}_\sigma^{(n)})}{n}=\infty?\]
Equivalently, is the linear maximal-overlap condition (Definition \ref{def: Return-Time Admissibility}) automatic for all finite measurable generators, or does it exclude some generators? The finite-window case has a linear bound (Proposition \ref{clopen cover proposition}), but the arbitrary measurable case is not settled by that argument.
\end{question}

\begin{question}[Global Bernoulli mean and variance]
\label{question:measure bernoulli}
For a nondegenerate Bernoulli shift, do we have
\[\Lambda_{\mu_{\mathbf p}}(\sigma)
=\frac1{H(\mathbf p)} \qquad \text{and/ or} \qquad \Lambda_{V,\mu_{\mathbf p}}(\sigma)=0?\]
If either equality fails, is there a way to determine the global value or establish whether the supremum is infinite.
\end{question}

\begin{question}[Uniform bounds across all admissible generators]
\label{question:uniform admissibility}
Can the constants in the inequality $M_\mu(\mathcal{P}_f^{(n)})\leq nC_{\mathcal{P}}$ be controlled uniformly over admissible generators for natural classes of systems? More generally, what dynamical hypotheses ensure that $\Lambda_\mu(f)$ and $\Lambda_{V,\mu}(f)$ are finite?
\end{question}

\begin{question}[Existence of admissible generators]
\label{question:admissible existence}
Which probability-preserving systems possess finite
return-time admissible generators? In particular, what is
the relationship between the existence of an ordinary
finite generator and the additional maximal-overlap bound?
\end{question}

\begin{question}[Zero-entropy and countable atomic systems]
\label{question:zero entropy}
Does there exist a probability-preserving system
with a finite return-time admissible generator
such that $\Lambda_\mu(f)=\infty$?
In particular, can this occur for a zero-entropy
system or a countable atomic system?
What values can $\Lambda_{V,\mu}(f)$ attain
in such examples?
\end{question}

\begin{question}[Distinguishing systems beyond entropy]
\label{question:beyond entropy}
Can the global return-time invariants distinguish positive-entropy measure-preserving systems with equal Kolmogorov-Sinai entropy? Can the variance invariant distinguish systems whose mean invariant agrees?
\end{question}

\begin{question}[Topological and measure-theoretic relations]
\label{question:topological measure relations}
Let $f:X\to X$ be continuous on a compact metric space and
let $\mu$ be an invariant Borel probability measure.
Under what additional hypotheses can the measure-theoretic
return-time invariants be compared with the topological
periodic-penetration invariants? Can either family distinguish
systems that the other does not?
\end{question}

\begin{question}[A finite generator with superlinear maximal overlap]
\label{question:generator superlinear}
Can one explicitly construct a finite measurable generator
of a two-sided Bernoulli shift for which the maximum
positive-measure overlap time of its forward joins grows
superlinearly? Any such construction must prove both the
generating property and the asserted zero-overlap conditions
for a sequence of join atoms. Small perturbations of a
canonical generator and approximation by finite-window
partitions alone do not settle the question.
\end{question}

\begin{question}[Classical recurrence and extremal overlap]
\label{question:overlap classical}
Can $\tau_\mu(B)$ and its distribution across generating
partitions be controlled by mixing coefficients, hitting-time
laws, return-time distributions, or quantitative recurrence
theorems? In particular, which hypotheses yield bounds that
are uniform over all atoms of a join, rather than merely
for each fixed measurable set?
\end{question}

\begin{question}[Alternative admissibility criteria]
\label{question:alternative-admissibility}
How do the resulting classes of generators and global invariants change if the maximal overlap condition $\displaystyle M_\mu(\mathcal P_f^{(n)})=O(n)$ is replaced by the weaker condition $\displaystyle A_\mu(\mathcal P_f^{(n)})=O(n)$? And, which admissibility conditions ensure meaningful and finite global invariants?
\end{question}

\bibliographystyle{plain}
\bibliography{references}

\vspace{1cm}

\noindent
\textbf{Diyath Nelaka Pannipitiya}\\
Department of Mathematical Sciences\\
Indiana University Indianapolis\\
Indianapolis, IN, USA\\
\textit{Email:}
\href{mailto:dinepann@iu.edu}{dinepann@iu.edu}

\end{document}